\documentclass[graybox]{svmult}

\usepackage{type1cm}        
\usepackage{makeidx}         
\usepackage{graphicx}        
\usepackage{multicol}        
\usepackage[bottom]{footmisc}

\usepackage{newtxtext}       %
\usepackage[varvw]{newtxmath}       
\usepackage{enumitem}

\usepackage{tikz}

\usepackage{mathtools}
\usepackage[colorlinks=true, allcolors=black]{hyperref}

\usepackage{subcaption}

\usepackage{multirow}

\newcommand{\cL}{\ensuremath{\mathcal{L}}}
\newcommand{\cM}{\ensuremath{\mathcal{M}}}
\newcommand{\cN}{\ensuremath{\mathcal{N}}}

\newcommand{\cP}{\ensuremath{\mathcal{P}}}

\newcommand{\cR}{\ensuremath{\mathcal{R}}}

\newcommand{\cV}{\ensuremath{\mathcal{V}}}

\newcommand{\bN}{\ensuremath{\mathbb{N}}}

\newcommand{\bR}{\ensuremath{\mathbb{R}}}

\newcommand{\rL}{\ensuremath{\mathrm{L}}}

\newcommand{\rU}{\ensuremath{\mathrm{U}}}
\newcommand{\rV}{\ensuremath{\mathrm{V}}}

\newcommand{\rY}{\ensuremath{\mathrm{Y}}}

\newcommand{\dx}{\ensuremath{\mathrm dx}}

\newcommand{\eps}{\ensuremath{\varepsilon}}
\newcommand{\dist}{\text{d}_V}
\newcommand*{\bary}{\text{Bar}}
\newcommand*{\Valpha}{\ensuremath{V^{(\alpha)}}}
\newcommand*{\cond}{\ensuremath{\; : \;}}
\newcommand*{\vspan}{\ensuremath{\text{span}}}
\newcommand*{\argmax}{\ensuremath{\text{arg}\, \text{max}}}
\newcommand*{\argmin}{\ensuremath{\text{arg}\, \text{min}}}
\renewcommand*{\vspan}{\ensuremath{\text{span}}}
\newcommand{\cVnlin}{\ensuremath{\cV_n^{\text{lin}}}}
\newcommand{\Vnlin}{\ensuremath{V_n^{\text{lin}}}}

\newcommand{\cMalpha}{\ensuremath{\cM^{(\alpha)}}}

\newcommand{\rYalpha}{\ensuremath{\rY^{(\alpha)}}}

\newcommand{\cVnbary}{\ensuremath{\cV_n^{\bary}}}
\newcommand{\rUalpha}{\ensuremath{\rU^{(\alpha)}}}
\newcommand{\greedy}{\ensuremath{\texttt{greedy}}}
\newcommand{\restartedgreedy}{\ensuremath{\texttt{restarted\_greedy}}}
\newcommand{\ealpha}{\ensuremath{e^{(\alpha)}}}

\renewcommand{\Pr}{\ensuremath{\cP_2(\Omega)}}
\newcommand{\ma}{\mu} 
\newcommand{\mb}{\nu}  
\newcommand{\mgI}{\nu} 

\newcommand{\barset}{\rU}
\newcommand{\barfun}{u}

\newcommand{\weight}{\Lambda} 
\newcommand{\weightcomp}{\lambda} 

\newcommand{\spt}{\ensuremath{\mathrm{spt}}}

\newcommand{\algoPG}{\ensuremath{\texttt{plain-greedy}}}
\newcommand{\algoY}{\ensuremath{$\rY$\texttt{-cart-split}}}
\newcommand{\algoM}{\ensuremath{$\cM$\texttt{-based-split}}}

\newcommand{\diffone}{\text{Diff-1}}
\newcommand{\difftwo}{\text{Diff-2}}
\newcommand{\convdiff}{\text{Cv-Diff}}
\newcommand{\kdv}{\text{KdV}}

\makeindex             

\begin{document}

\title*{Tree-Based Nonlinear Reduced Modeling}
\author{Diane Guignard and Olga Mula}
\institute{Diane Guignard \at Department of Mathematics and Statistics, University of Ottawa, Ottawa, ON K1N 6N5, Canada. 
\and Olga Mula \at Eindhoven University of Technology, Den Dolech 2, P.O. Box 513, 5600 Eindhoven, NL.}
\maketitle

\abstract{
This paper is concerned with model order reduction of parametric Partial Differential Equations (PDEs) using tree-based library approximations.
Classical approaches are formulated for PDEs on Hilbert spaces and involve one single linear space to approximate the set of PDE solutions. Here, we develop reduced models relying on a collection of linear or nonlinear approximation spaces called a library, and which can also be formulated on general metric spaces. To build the spaces of the library, we rely on greedy algorithms involving different splitting strategies which lead to a hierarchical tree-based representation. We illustrate through numerical examples that the proposed strategies have a much wider range of applicability in terms of the parametric PDEs that can successfully be addressed. While the classical approach is very efficient for elliptic problems with strong coercivity, we show that the tree-based library approaches can deal with diffusion problems with weak coercivity, convection-diffusion problems, and with transport-dominated PDEs posed on general metric spaces such as the $L^2$-Wasserstein space.}

\section{Introduction}

\subsection{Context}
This article develops tree-based nonlinear model order reduction of parametric Partial Differential Equations (PDEs). Let us start by drawing the general landscape. Parametric PDEs are of common use to model complex physical systems, and are routinely involved in design and decision-making processes.
Such equations can generally be written in abstract form as
\begin{equation}
\cP(u, y)=0,
\label{eq:ppde}
\end{equation}
where $\cP$ is a partial differential operator, and $y=(y_1,\dots,y_p)$ is a vector of
scalar parameters ranging in some domain $\rY\subset \bR^p$. The parameters represent physical quantities such as 
diffusivity, viscosity, velocity, source terms, or the geometry of the physical domain in which the PDE is posed.  In several relevant instances, $y$ may be high or even countably infinite dimensional, that is, $p\gg1$ or $p=\infty$ (for simplicity, we consider $p<\infty$ in this work).

We assume well-posedness, that is, for any $y\in \rY$  the problem admits a unique solution $u=u(y)$ in a metric space $(V, \dist)$ whose elements depend on a physical variable $x$ ranging in a domain $\Omega\subset \bR^q$.
The variable $x$ usually refers to space but it is not limited to that meaning, and it may also refer to more elaborate sets of variables such as space and time. We may thus regard $u$ as a function $(x,y)\mapsto u(x, y)$ from $\Omega\times \rY$ to $\bR$,
or we may also consider the \emph{parameter-to-solution map}
\begin{equation}
y\in \rY\mapsto u(y)\in V
\label{eq:solmap}
\end{equation}
from $\rY$ to $V$. This map is typically nonlinear, and the image set 
\begin{equation}
\cM\coloneqq  u (\rY) =\{u(y) \, : \, y\in Y\}\subset V
\label{eq:manifold}
\end{equation}
is usually called  the {\it solution manifold}. Throughout the paper, we assume that $\rY$ is compact in $\bR^p$ and that the map \eqref{eq:solmap}
is continuous. Therefore, $\cM$ is a compact set of $V$. 

In numerous design and decision-making tasks, one is often confronted to optimize over the solution manifold $\cM$, and this requires evaluating solutions $u(y)$ on a large number of dynamically updated parameters $y\in \rY$. Computations cannot be addressed rapidly unless the overall complexity has been appropriately reduced. This motivates the search for accurate methods to approximate the family of solutions very quickly at a reduced computational cost. This task, usually known as \emph{reduced modelling} or \emph{model order reduction}, has classically been considered in Hilbert and Banach spaces $(V, \Vert \cdot\Vert_V)$ by approximating $\cM$ with well-chosen linear subspaces of $V$. In other words, one finds an appropriate $n$-dimensional subspace $\Vnlin\subset V$, and builds a map
$$
\hat u_n: \rY \to \Vnlin, \quad y \mapsto \hat u_n(y),
$$
in which evaluations of $\hat u_n$ scale only with $n$ in complexity. This approach is suitable if the approximation error
$$
e(\cM, \hat u_n, V) \coloneqq \sup_{y\in \rY} \dist(u(y), \hat u_n(y))
$$
decays fast as $n\to \infty$. This is typically the case if $\hat u_n$ has a near-optimal approximation property in the sense that there exists $C(n)>0$ such that
$$
e(\cM, \hat u_n, V) \leq C(n) \dist(\cM, \Vnlin),
$$
where
\begin{equation}
\label{eq:err-best-approx-manifold}
\dist(\cM, \Vnlin) \coloneqq \sup_{u\in \cM} \inf_{v\in \Vnlin} \dist(u, v)
\end{equation}
is the error of best approximation of $\cM$ in $\Vnlin$. If the product $C(n)\dist(\cM, \Vnlin)$ decays fast as $n\to \infty$, the approximation of the parameter-to-solution map \eqref{eq:solmap} by $\hat u_n$ is computationally very efficient. However, we can expect this to happen only when $C(n)$ mildly increases with $n$, and the Kolmogorov $n$-width of $\cM$
\begin{equation}
\label{eq:kolmo-width}
\kappa(\cM, \cVnlin, V) \coloneqq \inf_{\Vnlin \in \cVnlin} 
\dist(\cM, \Vnlin)
\end{equation}
decays fast with $n$. In the above definition, $V_n^{\text{lin}}$ runs over the set $\cVnlin$ of all linear subspaces of dimension at most $n$ in the Banach space $V$. While this is the case for certain families of parabolic or elliptic problems (see \cite{CD2015}), most transport-dominated problems are expected to present a slow decaying width and require to study nonlinear approximation methods (see, e.g., \cite[Chapter 3]{BCOW2017}).

\subsection{Library approximation}
There is a tremendous activity in overcoming the limitations of working with linear subspaces for model reduction, and this paper is a contribution in this direction. We formulate a library approach which can be expressed in general metric spaces, and for which one can recover more classical formulations when instantiating it for PDEs on Hilbert and Banach spaces.

The idea is as follows. Given a metric space $V$, one considers an approximation family $\cV_n$ which is not necessarily restricted to linear subspaces $\cVnlin$ as in the above discussion. Each member $V_n\in \cV_n$ can in general denote a linear or nonlinear approximation space parametrized by at most $n$ degrees of freedom, and we define
\begin{equation}
\dist(\cM, V_n)\coloneqq \sup_{u\in \cM} \inf_{v\in V_n} \dist(u, v)
\label{eq:err-best-approx-manifold-nonlinearVn}
\end{equation}
and
$$
\kappa(\cM, \cV_n, V)\coloneqq \inf_{V_n \in \cV_n} 
\dist(\cM, V_n)
$$
for such general spaces $V_n$ exactly in the same spirit as \eqref{eq:err-best-approx-manifold} and \eqref{eq:kolmo-width}. The quantity $\kappa(\cM, \cV_n, V)$ becomes a generalized Kolmogorov $n$-width which boils down to the classical definition when $V$ is a Banach space and $\cVnlin$ is the family of $n$-dimensional linear subspaces. Among the many possible examples of nonlinear $\cV_n$, we may mention a family $\cV_n^{\cN\cN}$ generated by neural network functions with different architectures. Each $V_n \in\cV_n^{\cN\cN} $ would be the set of neural network functions with a given architecture and having at most $n$ parameters varying in a certain range. Here, the variability of the architectures generates $\cV_n$ in a Hilbert or a Banach space. Another example, which we will develop later on, is when $V$ is the $L^2$-Wasserstein space of measures. In this case, we can consider $\cV_n=\cVnbary$ as the set of all possible barycenters generated by at most $n$ measures. This set will be defined more precisely in equation \eqref{eq:cVnbary}.

We next go beyond this idea with a library approximation. Instead of finding only one space $V_n\in \cV_n$ for approximating $\cM$, we consider a collection of spaces, which we call a library $\cL = \{ \Valpha\}_{\alpha\in A}$ of cardinality $m\coloneqq|\cL|=|A|\in \bN^*$, and for which each $\Valpha$ is a linear or nonlinear space of dimension $n(\alpha)$. In our work, we rely on greedy algorithms to build spaces with different splitting strategies which will eventually lead to a representation of $\cL$ in terms of a tree of subspaces. This is interesting because the representation of the subspaces $\Valpha\in \cL$ is thus hierarchical. The goal is then to explore how much we can benefit from this compression depending on different tree constructions.

To keep track of the tree structure, $\alpha$ will denote a multi-index set of the form $\alpha = (a_1,\dots, a_k)\in \bN_0^k$ for some $k\in \bN^*$, and $A$ is the set of multi-indices that are needed to define $\cL$. Once $\cL$ is obtained, we build an approximation
\begin{equation}
\label{eq:library-approach}
\hat u_\cL : \rY \to \cL, \quad \hat u_\cL(y)\in V^{(\iota(y))}, 
\end{equation}
where $\iota : \rY \to A$ is a mapping from the parameters to the subspace index.

In the same spirit as the notions $\dist(\cM, V_n)$ and $\kappa(\cM, \cV_n, V)$, we can now define the error of best approximation of $\cM$ with a given library $\cL$ as
$$
\dist(\cM, \cL) \coloneqq \sup_{u\in \cM} \inf_{\Valpha\in \cL} \inf_{v\in \Valpha} \dist(u, v)
$$
and a notion of library width as
\begin{equation}
\label{eq:library-width}
\kappa(\cM, \cL_{m, n}, V) \coloneqq \inf_{\cL \in \cL_{m,n}} 
\dist(\cM, \cL)
\end{equation}
with
\begin{equation}
\label{eq:library-set}
\cL_{m,n} \coloneqq \{ \cL = (\Valpha)_{\alpha\in A} \cond \vert \cL \vert \leq m, \text{ and } n(\alpha)\leq n \text{ for all }\alpha \in A \}
\end{equation}
the set of libraries of cardinality at most $m$ and containing subspaces $\Valpha$ of dimension $n(\alpha)\leq n$. This approach is expected to give better accuracy than the one without a library since we have 
\begin{equation} \label{rel:widths}
	\kappa(\cM, \cL_{m, n}, V)
	\leq \kappa(\cM, \cL_{1, n}, V) = \kappa(\cM, \cV_n, V).
\end{equation}

\subsection{Contribution and connection with earlier works}
The idea of library approximation is not new in model order reduction. However, the common denominator of most approaches is that they rely on the assumption that the parameter space $\rY$ is a cartesian product or has a very simple shape (we refer to \cite{EPR2010,MS2013,BCDGJP2021,GGLM2021, CDMN2022, DCAKR2022, NP2023, GJ2024} for works in this direction). The construction of the library then relies on identifying appropriate splittings of $\cM$ induced by rectangular splittings of $\rY$ (or splittings involving a very explicit expression in $\rY$). While this type of splitting is very practical because the definition of the index mapping $\iota$ is straightforward, this point of view cannot be easily extended for general $\rY$, and there is no real exploitation of a potential tree structure in the representation of the library. With respect to this body of works, one of the main novelties in this paper is therefore the introduction of a splitting algorithm that is valid for any shape of the parameter domain and the compressed tree representation. The algorithm that we explore actually emerged in a discussion between Prof.~Ron DeVore and Diane Guignard (one the authors of this paper).

The second main contribution is the phrasing of our algorithms in the general setting of nonlinear metric spaces. While the lack of linearity of $V$ may pose numerical challenges, working in nonlinear metric spaces comes with the significant advantage of potentially capturing physical effects intrinsically in the geometry of the space. This gives room to reduce complexity even further, preserve certain physical properties by construction (e.g.~mass conservation, symplecticity,...), and to treat different PDE classes under the same functional framework. A salient example is the family of Wasserstein gradient flows for which it is proven that some of the most common PDEs such as the heat equation or the Fokker-Planck equation can be seen as gradient flows in the Wasserstein space $(\cP_2(\Omega), W_2)$ (which we will define more rigorously later on).

Probably \cite{ELMV2020} is one of the first works exploring model reduction on metric spaces where a strategy for the $L^2-$Wasserstein space was introduced. In that work, a parametric family of Wasserstein barycenters plays the role of $V_n$, and one can define a barycentric $n$-width exactly in the same manner as \eqref{eq:kolmo-width} but where $V_n$ is infimized over all barycenter approximation sets $\cV_n$ generated by at most $n$ measures. We refer to \cite{BBEELM2023, DDEL2023, Blickhan2023, DFM2023} for further works and applications of the idea. The work \cite{DFM2023} is of particular interest for this paper since it develops a fully adaptive library approach. The main difference with the present approach is that it is not tree-based, and it involves a much larger library $\cV_n$. This makes the definition of the operator $\iota$ much more costly. We discuss more in detail the differences between both approaches later on in Section \ref{sec:barycentric-greedy}.

\subsection{Plan of the paper}
The paper is organized as follows. Since our tree-based library construction heavily relies on greedy algorithms, we recall in Section \ref{sec:plain-greedy} the vanilla versions of this algorithm for $\cM$ in Banach spaces, and in the $L^2$-Wasserstein space. In Section \ref{sec:tree-based-library-greedy} we explain our tree-based library approaches. Section \ref{sec:numerics} contains numerical experiments illustrating the performances of the algorithms. 

\section{The plain greedy algorithm for model reduction}
\label{sec:plain-greedy}
In this section, we recall the greedy algorithms for model reduction that we use for our library approximation. The first is the classical greedy algorithm in Hilbert and Banach spaces whose convergence properties were analyzed in \cite{BCDDPW2011, DPW2013}. The second greedy algorithm is less standard, and was originally introduced in \cite{ELMV2020} for model reduction in the $L^2$-Wasserstein space.

\subsection{Formulation in Hilbert and Banach spaces}
\label{sec:vanilla-greedy}
Let $(V, \Vert \cdot\Vert_V)$ be a Banach space, and consider $\cV_n=\cVnlin$. The greedy algorithm builds a space $\Vnlin\in \cVnlin$ in a nested, iterative way. For $n=1$, we find
\begin{equation}
\label{eq:vanilla-greedy-n1}
u_1 \in \argmax_{u\in \cM} \Vert u \Vert_V
\end{equation}
and define
$$
\rU_1 \coloneqq \{u_1\}, \quad V_1^{\text{lin}}\coloneqq \vspan\{u_1\}.
$$
For $n>1$, given $\rU_{n-1}$ and $V_{n-1}^\text{lin}$, we find
\begin{equation}
\label{eq:vanilla-greedy-n}
u_n \in \argmax_{u\in \cM} \min_{v\in V_{n-1}^\text{lin}} \Vert u -v\Vert_V
\end{equation}
and set
$$
\rU_n = \rU_{n-1} \cup \{u_n\}, \quad \Vnlin = \vspan\{ \rU_n \}.
$$
In \cite{BCDDPW2011, DPW2013}, it was proven that this greedy procedure generates a sequence of spaces $\{\Vnlin\}_n$ for which $\dist(\cM, \Vnlin)$ decays at a comparable rate as the Kolmogorov $n$-width $\kappa(\cM, \cV_n^\text{lin}, V)$ as $n\to \infty$.

Before going further, several comments on the practical implementation of the algorithm are in order. Note that we have assumed that the above optimization problems over $\cM$ admit a maximum. In general, this is not true and, in addition, one may not have a perfect optimization procedure to find the exact maximum. We can overcome both issues by formulating a weak greedy version of \eqref{eq:vanilla-greedy-n1}-\eqref{eq:vanilla-greedy-n} where, for a parameter $\gamma\in (0, 1]$, we search for
\begin{equation}
u_1 \in \{ u \in \cM\cond \Vert u \Vert_V \geq \gamma \sup_{v\in \cM} \Vert v \Vert_V  \}
\end{equation}
and for $n>1$,
\begin{equation}
\label{eq:weak_greedy_err}
u_n \in \{ u \in \cM\cond \inf_{v\in V_{n-1}^\text{lin}}\Vert u-v \Vert_V \geq \gamma \sup_{z\in \cM} \inf_{v\in V_{n-1}^\text{lin}} \Vert z-v \Vert_V  \}.
\end{equation}
When $\gamma=1$, we recover the original greedy algorithm.

The implementation of the weak greedy algorithm is still computationally expensive in its current form because it requires having access to all $u\in \cM$, and this means that we need to compute the solution $u(y)$ for all parameters $y\in \rY$. The complexity can be greatly reduced when working with PDEs defined in Hilbert spaces via weak formulations. In this case, we can replace the computation of full PDE solutions and distance evaluations by evaluation of residuals. The justification of this is as follows. In the Hilbert space setting, \eqref{eq:ppde} can be expressed as an operator equation of the form $\cP(u, y)=A(y)u(y)-f(y)=0$. Here the solution $u(y)$ lives in the solution Hilbert space $V$, we have identified a test Hilbert space $Z$ such that $f(y)$ belongs to the dual $Z'$ of $Z$, and the operators are boundedly invertible mappings from $V$ to the dual $Z'$. We assume that the invertibility is uniform in $y\in \rY$, that is, there exists $0<r\leq R<\infty$ such that
\begin{equation}
\label{eq:boundedinv}
R(y)\coloneqq \Vert A(y)\Vert_{V\to Z'} \leq R, \quad \text{and}\quad
r(y)\coloneqq \Vert A(y)^{-1}\Vert_{Z'\to V} \leq r^{-1}, \quad \forall y\in \rY.
\end{equation}
This setting covers a wide range of PDEs, such as classical elliptic problems with $Z=V$, as well as saddle-point problems, indefinite problems and unsymmetric problems such as convection-diffusion PDEs, or space-time formulations of parabolic problems. For some of these PDEs, $Z$ sometimes needs to be chosen different from $V$ to ensure a moderate value of the condition number $\nu \coloneqq R/r\geq 1$. For example, for a simple elliptic operator of the form $A(y)u(y)=-\text{div}(a(y)\nabla u(y))$ where the parametric diffusion $a(y)$ is uniformly bounded away from $0$ and $\infty$, then the constants $r$ and $R$ directly relate to $\min_{y\in \rY} a(y)$ and $\max_{y\in \rY} a(y)$. It follows from \eqref{eq:boundedinv} that one has the equivalence
\begin{equation} \label{eq:equiv_res_err}
r\Vert u(y) - v \Vert_V
\leq \cR(v, y) \leq R\Vert u(y) - v \Vert_V,	
\end{equation}
where
$$
\cR(v,y)\coloneqq \Vert A(y) v -f(y)\Vert_{Z'}
$$
is the residual of the PDE for a state $v\in V$ and a parameter $y\in \rY$. As a consequence, one can further modify the weak greedy algorithm into the search for $n>1$ of
\begin{equation}
	\label{eq:weak_greedy_est}
	y_n \in \{ y\in \rY \cond \Delta_{n-1}(y) \geq \gamma \sup_{y\in \rY} \Delta_{n-1}(y)  \},
\end{equation}
where
$$
\Delta_{n-1}(y)\coloneqq\inf_{v\in V_{n-1}^\text{lin}}\cR(v, y)
$$
is a quantity that one can estimate without computing exact solutions of the PDE.

Once $y_n$ is obtained from \eqref{eq:weak_greedy_est}, we set $u_n = u(y_n)$. Thanks to the equivalence relation \eqref{eq:equiv_res_err}, $u_n$ satisfies a weak greedy property with respect to the ambient error norm,
$$
\inf_{v\in V_{n-1}^\text{lin}}\Vert u_n-v \Vert_V \geq \widetilde \gamma \sup_{z\in \cM} \inf_{v\in V_{n-1}^\text{lin}} \Vert z-v \Vert_V,
$$
for a parameter
$$
\widetilde \gamma \coloneqq \frac{\gamma}{\nu} \leq \gamma.
$$
This means that working with the equivalent residual degrades the tightness of the weak greedy algorithm by a factor $\nu$.

Step \eqref{eq:weak_greedy_est} only requires computing residuals that do not directly involve PDE solutions. However, we still have to face the fact that we need to make a search over the whole parameter set $\rY$. Depending on the type of residual, one directly deals with the optimization task with classical optimization methods. This is the case  when the residual depends affinely on the parameters. For more complex dependencies, the search over $\rY$ is replaced by a search over a large, finite training set $\rY^{(N)}$ of parameters (which generates a sampled version $\cM^{(N)}$ of the manifold). This is a change that affects the global approximation properties of the algorithm, and quantifying how many samples are enough to preserve the same convergence rate as the non-sampled version is a difficult question. Some elements of answer can be found in \cite{CDDN2020}.

\begin{remark} \label{rem:y_dep}
When an explicit formula for the $y$-dependent constant $r(y)$ in \eqref{eq:boundedinv} is available, then it is customary to use the \emph{a posteriori} error estimator $\Delta_n(y)/r(y)$ rather than $\Delta_n(y)$ in the weak greedy procedure, see for instance \cite{QMN2016,HRS2016}. Indeed, we have
$$\frac{1}{R(y)}\Delta_n(y)\le e_n(y) \le \frac{1}{r(y)}\Delta_n(y), \quad e_n(y)\coloneqq\inf_{v\in V_n^\text{lin}}\|u(y)-v\|_V,$$
and thus this estimator provides a certified control of the error with constant one, and its effectivity index belongs to $[1,R(y)/r(y)]$ (the effectivity index of $\Delta_n(y)$ lies in $[r(y),R(y)]$). In practice, including the constant $r(y)$ or not can impact the selection of the snapshots, in particular when $r(y)$ varies significantly with $y$. However, this sensitivity will not be investigated in this work.
\end{remark}

For subsequent developments, for a given Banach space $V$, a given target accuracy $\eps>0$ and a given maximal amount of iterations $n_{\max}$, it will be convenient to define the routine
\begin{equation}
\label{eq:greedy-V-Banach}
[\rU_n, \Vnlin, e] = \greedy[\cM, V, \eps, n_{\max}]
\end{equation}
which iterates over the steps \eqref{eq:vanilla-greedy-n1}-\eqref{eq:vanilla-greedy-n} of the greedy algorithm until
$$
n = \min\{ n_{\max}, n_{\eps} \}
$$
with
$$
n_\eps = \min\{ n\in \bN^* \cond \dist(\cM, V_n)\leq \eps \}.
$$
This routine yields the space $\Vnlin=\vspan\{\rU_n\}$ which approximates $\cM$ with accuracy $e=\dist(\cM, V_n)\leq \eps$ provided that reaching this accuracy takes less than $n_{\max}$ steps. If it takes more than $n_{\max}$ iterations to reach $\eps$ accuracy, then $e=\dist(\cM, V_n)>\eps$ and the approximation error does not meet the target accuracy.

It will also be necessary to introduce a restarted version of the greedy algorithm where we resume iterations from a given space $V_k^{\text{lin}}$ with $k\leq n_{\max}$, and we compute the iterations \eqref{eq:vanilla-greedy-n} starting at iteration $k+1$. This routine is denoted
\begin{equation}
\label{eq:restarted-greedy-V-Banach}
[\rU_n, \Vnlin, e] = \restartedgreedy[\cM, V, V_k^{\text{lin}}, \eps, n_{\max}].
\end{equation}
and the output is of dimension $n \in [k, n_{\max}]$.

\subsection{A greedy barycentric algorithm in the $L^2$-Wasserstein space}
\label{sec:vanilla-barycentric-greedy}
The above weak greedy algorithm can be formulated in an analogous manner in the $L^2$-Wasserstein space of measures using barycenters to approximate measures. To explain this, we first recall the definition of this space and of barycenters.

\subsubsection{The $L^2$-Wasserstein space and barycenters}
\label{sec:wasserstein-space}

We assume in the following that $\Omega$ is a subset in the normed vector space $(\bR^q, \vert \cdot \vert)$, with $ \vert \cdot \vert$ the Euclidean norm on $\bR^q$. We denote $\cM(\Omega)$ the space of Borel regular measures on $\Omega$ with finite total mass and
\begin{equation*}
\cM^+(\Omega)\coloneqq\{ \mgI \in \cM(\Omega) \cond \mgI \geq 0 \}~,\quad
\cP(\Omega)\coloneqq\{ \mgI \in \cM^+(\Omega) \cond \mgI(\Omega)=1 \}~.
\end{equation*}
The Wasserstein space $\cP_2(\Omega)$ is defined as the set of probability measures $\mgI\in \cP(\Omega)$ with finite second order moments, namely
$$
\Pr ~\coloneqq~ \{ \mgI \in \cP(\Omega) \cond \int_\Omega \vert x \vert^2 \,\mgI(\dx) \;< +\infty \}~.
$$
We endow the output space $\cP_2(\Omega)$ with the metric that is induced by optimal transport theory which is defined as follows. Let $\ma$ and $\mb$ be two probability measures in $\cP_2(\Omega)$. We define $\Pi(\ma, \mb)\subset \cP_2(\Omega\times \Omega)$ as the subset of probability distributions $\pi$ on $\Omega\times \Omega$ with marginal distributions equal to $\ma$ and $\mb$. The Wasserstein distance between $\ma$ and $\mb$ is defined as:
\begin{equation}
\label{eq:W2}
W_2(\ma,\mb) ~\coloneqq~ \mathop{\inf}_{\pi \in \Pi(\ma,\mb)} \left( \int_{\Omega \times \Omega} \vert x-y\vert^2 \,\pi(dx,dy) \right)^{1/2},\quad \forall (\ma,\mb)\in \cP_2(\Omega) \times \cP_2(\Omega)~.
\end{equation}
As detailed in \cite{Villani2003}, the space $\cP_2(\Omega)$ endowed with the distance $W_2$ is a metric space that is usually called the $L^2$-Wasserstein space.

In this space, Wasserstein barycenters arise as natural objects to define a reasonable, and computationally tractable approximation family $\cV_n$ in $V=(\Pr, W_2)$. Their definition and properties are relatively well understood, especially since the works of \cite{AC2011}. Let $n\in \bN^*$ denote the number of observations in our dataset and consider the simplex in $\bR^n$:
$$
\Sigma_n ~\coloneqq~ \Big\{\, \weight_n = (\weightcomp_1,\dots, \weightcomp_n)^T \in \bR^n\cond \weightcomp_i\geq 0,\, \sum_{i=1}^n \weightcomp_i = 1 \, \Big\}~.
$$
Given a set of weights $\weight_n = (\weightcomp_i)_{1\leq i\leq n} \in \Sigma_n$, and given a set $\barset_n = \{\barfun_i\}_{1\leq i\leq n}$ of $n$ probability measures from $\cP_2(\Omega)$, we say that $\bary(\weight_n, \barset_n) \in \cP_2(\Omega)$ is a barycenter associated to $\weight_n$ and $\barset_n$ if and only if
\begin{equation}
\label{eq:barygen}
\bary(\weight_n, \barset_n)
\;\in\;
\argmin_{\mgI \in \cP_2(\Omega)} \sum_{i=1}^n \weightcomp_i W_2^2(\mgI,\barfun_i)~.
\end{equation}
We can thus define an approximation space $V_n^\bary$ in $\Pr$ by varying the weights $\Lambda_n\in \Sigma_n$,
\begin{equation}
V_n^\bary \coloneqq \bary(\Sigma_n, \rU_n)
~\coloneqq~ \{  \bary(\weight_n, \barset_n) \cond \weight_n \in \Sigma_n \} \subset \cP_2(\Omega)~,
\label{eq:bary}
\end{equation}
and this space belongs to the approximation class
\begin{equation}
\label{eq:cVnbary}
\cV_n^{\bary} \coloneqq \{ \bary(\Sigma_n, \rU_n) \cond \rU_n \in \Pr^n \}.
\end{equation}

\subsubsection{The barycentric greedy algorithm}
\label{sec:barycentric-greedy}
If $\cM \subset V=(\Pr, W_2)$, we can build a barycentric space $V_n=\bary(\Sigma_n, \rU_n)$ to approximate $\cM$ by building $\rU_n$ in a nested manner with a greedy algorithm. This is done exactly in the same spirit as the construction for the Hilbert/Banach case.

The algorithm is initialized with the two elements of $\cM$ that lie at furthest distance, namely
$$
(u_1, u_2) \in \argmax_{(\nu_1,\nu_2)\in \cM\times\cM} W_2(\nu_1,\nu_2),
$$
and we set
$$
\rU_2 \coloneqq \{u_1,u_2\}, \quad V_2^\bary \coloneqq \bary(\Sigma_2, \rU_2).
$$
Next, for $n\geq 3$, given $\rU_{n-1}$ and the associated $V_{n-1}^\bary=\bary(\Sigma_{n-1},\rU_{n-1})$, we search for
$$
u_n \in \argmax_{\nu \in \cM} \inf_{\zeta\in V_{n-1}^\bary} W_2(\nu, \zeta),
$$
and we set
\begin{equation}
\label{eq:barycentric-Vn}
\rU_n=\rU_{n-1}\cup\{u_n\}, \quad V_n^\bary \coloneqq \bary(\Sigma_n, \rU_n).
\end{equation}

The convergence properties of this algorithm are mathematically not well understood yet but there is ample numerical evidence that, for many PDE classes, the algorithm builds spaces $V_n^\bary$ that efficiently approximate $\cM$ in the sense of presenting a fast convergence rate $\dist(\cM, V_n^\bary)$ (see equation \eqref{eq:err-best-approx-manifold-nonlinearVn}). Like in the Banach space case, the practical implementation of the algorithm requires replacing $\cM$ by a finite set of $N\gg1$ training samples $\rU^{(N)} =\{u(y) \cond y \in \rY^{(N)}\}$ computed from a finite set of $N$ parameters $\rY^{(N)}$. Computations present an additional challenge compared to the Banach space setting which is connected to the complexity of evaluating Wasserstein distances and barycenters. In one space dimension ($q=1$), one can leverage closed form characterizations involving inverse cumulative distributions. This is no longer possible for $q>1$ and, in this setting, a popular approach is to resort to entropy regularized versions of the distance and the barycenter (see \cite{PC2019, DFM2023}).

We can now easily adapt the definition of the $\greedy$ and $\restartedgreedy$ routines introduced in \eqref{eq:greedy-V-Banach} and \eqref{eq:restarted-greedy-V-Banach} for $V=(\Pr, W_2)$. We just need to replace $\Vnlin$ by the barycentric space $V_n^\bary$ from \eqref{eq:barycentric-Vn}.

A library version of this algorithm was recently introduced in \cite{DFM2023}. We briefly recall it next in order to compare it to the tree-based library approach that we introduce in the next section. The starting point in that work is to observe that, since in practice we use a training set $\rU^{(N)}$ for the greedy algorithm, the $V_n$ that is obtained in the plain greedy algorithm belongs to the approximation class
\begin{equation}
\label{eq:n-sparse-bary}
\cV^{\bary}_{n, N} ~\coloneqq~ \bary(\Sigma_N^n, \rU^{(N)}) ~=~ \{ \bary(\weight_N, \rU^{(N)}) \cond \weight_N \in \Sigma_N^n \}
\subset \cV_n^\bary~,
\end{equation}
where
$$
\Sigma_N^n ~\coloneqq~ \{ \weight_N\in \Sigma_N  \cond \vert \spt(\weight_N)\vert \leq n \}
$$
is the set of $n$-sparse vectors from the simplex $\Sigma_N$. Here, $\vert \spt(\weight_N)\vert$ denotes the cardinality of the support of $\weight_N$, namely
$$
\spt(\weight_N) ~\coloneqq~ \{i \in \{1,\dots, N\} \cond \weightcomp_i \neq 0\}~.
$$
It is interesting to remark that $\Sigma_N^n$ is the union of the $n$-simplices embedded in $\bR^N$, and there are $\binom{N}{n}$ of these simplices. This means that $|\cV^{\bary}_{n, N}|= \binom{N}{n}$. In \cite{DFM2023}, the authors build a model reduction method which can be summarized as the construction of a mapping
$$
\hat u_{\cV^{\bary}_{n, N}} : \rY \to \cV^{\bary}_{n, N}\,, \quad u_{\cV^{\bary}_{n, N}}\in V_n^\bary(y) \in \cV^{\bary}_{n, N}.
$$
The mapping involves a procedure to adapt the approximation space $V_n^\bary(y)$ for each $y\in \rY$. This is done via solving a challenging optimization problem involving sparsity constraints (we refer to \cite{DFM2023} for further details). As a consequence, this is a library approach of the form \eqref{eq:library-approach} where $\cL\coloneqq\cV^{\bary}_{n, N}$. One of the main differences compared to the strategy that we develop next is that, for each $y\in \rY$, the search of $V_n^\bary(y)$ is done over the whole $\binom{N}{n}$ members of the family which is very large if $N\gg1$. In this context, our tree-based approach can be understood as a way of narrowing-down the amount of spaces in which one needs to make the search.

\section{Tree-based library approaches}
\label{sec:tree-based-library-greedy}
We next present two tree-based library approaches that build on the previous greedy algorithms. The strategy of Section \ref{sec:cartesian-split} requires $\rY$ to have a tensor product structure. The second approach, which is one of the main novelties of the paper, can be applied to any shape of parameter domain $\rY$. We present it in Section \ref{sec:general-split}.

\subsection{Cartesian splittings of $\rY$}
\label{sec:cartesian-split}
One can build on the plain greedy algorithm to derive a library approach where $\cL=\{\Valpha\}_{\alpha\in A}$ is composed of spaces $\Valpha$ that are built in order to approximate a subset $\cMalpha$ of the manifold in such a way that
$$
\dist(\cMalpha, \Valpha)\leq \eps, \quad \forall \alpha\in A.
$$
The subsets $\cMalpha$ must form a covering of $\cM$ in the sense that
$$
\cM=\cup_{\alpha\in A} \cMalpha.
$$
Since $\cM$ has a natural parametrization through $\rY$, the most natural approach is to build a splitting of $\rY$,
$$
\rY=\cup_{\alpha\in A}\rYalpha,
$$
which we then use to build
$$
\cMalpha = \{ u(y)\cond y \in \rYalpha \}.
$$
Splitting $\rY$ is particularly easy when it is a cartesian product since we can build an approach based on simple dyadic partitions. To explain the algorithm, we assume without loss of generality that $\rY=[0, 1]^p$, and we define indices
\begin{equation} \label{eq:indices_alpha}
	\alpha = \{(k_m, i_m) \}_{m=1}^p, \quad \text{for }k_m\in \bN_0, \;0\leq i_m < 2^{k_m}.
\end{equation}
We associate to $\alpha$ the parameter subdomain 
$$
\rYalpha \coloneqq \times_{m=1}^p \left[ \frac{i_m}{2^{k_m}}, \frac{i_m+1}{2^{k_m}}\right]\; \subseteq \rY.
$$
Note that with this notation
$$
\rY = \rY^{(\alpha)}, \; \cM = \cMalpha\quad \text{ for } \alpha = \{(0, 0) \}_{m=1}^p,
$$
and for a given $\alpha= \{(k_m, i_m) \}_{m=1}^p$ we can express partitions of $\rY^{(\alpha)}$ along a certain parameter $\tilde p\in\{1,\ldots,p\}$ as
$$
\rY^{(\alpha)} = \rY^{(\alpha(0, \tilde p) )} \cup \rY^{(\alpha(1, \tilde p) )}
$$
with
\begin{equation*}
\begin{cases}
\alpha(0, \tilde p)
&= \{ (k_1, i_1), \dots, (k_{\tilde p -1}, i_{\tilde p-1}), (k_{\tilde p}+1, 2i_{\tilde p}),(k_{\tilde p +1}, i_{\tilde p+1}), \dots, (k_p, i_p)  \} \\
\alpha(1, \tilde p)&= \{ (k_1, i_1), \dots, (k_{\tilde p -1}, i_{\tilde p-1}), (k_{\tilde p}+1, 2i_{\tilde p}+1),(k_{\tilde p +1}, i_{\tilde p+1}), \dots, (k_p, i_p)  \}.
\end{cases}
\end{equation*}
Note that for $m\in\{1,2,\ldots,p\}$, $k_m$ controls the length of the subinterval in the direction $y_m$ while $i_m$ controls its location in $[0,1]$.

We can now formulate the following algorithm to generate a tree-based library approximation. We set the target tolerance $\eps>0$ and the maximal dimension $n_{\max}$ of the elements of the library. We then proceed iteratively as follows:
\begin{itemize}
\item Initialization ($k=1$): We define
\begin{align*}
\alpha_k &= \{(0, 0) \}_{m=1}^p \\
A_k &= \{\alpha_k\} \\
\rY_k &= \cup_{\alpha \in A_k} \rY^{(\alpha)} = \{ \rY \}
\end{align*}
and we run
$$
[\rU^{(\alpha_k)}, V^{(\alpha_k)}, e^{(\alpha_k)}]
= \greedy(\cM^{(\alpha_k)}, V, \eps, n_{\max}),
$$
where, by construction, $\cM^{(\alpha_k)}=\cM$. We then set
$$
\rL_k=\{ \Valpha\}
$$
which is a library containing one space of dimension at most $n_{\max}$. Therefore, $\rL_k \in \cL_{1, n_{\max}}$ according to Definition \eqref{eq:library-set}.

If $e^{(\alpha_k)}\leq \eps$, the algorithm stops and we set $K=1$. Otherwise we set $k=2$ and the algorithm goes into the inductive step.
\item Inductive step $k>1$: From the previous step, we receive $A_{k-1},\,\rY_{k-1}$, and $\rL_{k-1}=\{\Valpha\}_{\alpha\in A_{k-1}}$. We start by setting $A_{k}=A_{k-1}$, and we are going to update the definition of $A_k$ through the following procedure. For every $\alpha\in A_{k-1}$:
\begin{itemize}
    \item If $e^{(\alpha)}=\dist(\cMalpha, \Valpha)\leq\eps$, the algorithm has reached accuracy $\eps$ on $\cMalpha$ so we do not need to split this subset.
    \item If $e^{(\alpha)}=\dist(\cMalpha, \Valpha)>\eps$, we  split $\cMalpha$ along one parameter direction. To discover the best parameter to choose, we run a greedy algorithm for every possible dyadic split. This means that for every parameter $\tilde p \in \{1, \dots,p\}$, we compute
\begin{equation*}
\begin{cases}
[\rU^{(\alpha(0, \tilde p))}, V^{(\alpha(0, \tilde p))}, e^{(\alpha(0, \tilde p))}]
&= \greedy[\cM^{((\alpha(0, \tilde p))}, V, \eps, n_{\max}] \\
[\rU^{(\alpha(1, \tilde p))}, V^{(\alpha(1, \tilde p))}, e^{(\alpha(1, \tilde p))}]
&= \greedy[\cM^{(\alpha(1, \tilde p))}, V, \eps, n_{\max}]
\end{cases}
\end{equation*}
and we define
\begin{equation} \label{crit_psplit}
\ealpha_{\tilde p} \coloneqq \max( e^{(\alpha(0, \tilde p))}, e^{(\alpha(1, \tilde p))} ).	
\end{equation}
We next choose the 
parameter along which we do the splitting by finding
\begin{equation} \label{eq:p_split}
p_{\text{split}}
\coloneqq \argmin_{\tilde p\in \{1,\dots, p\}} \ealpha_{\tilde p}.	
\end{equation}
Then we decompose $\rY^{(\alpha)} \in \rY_{k-1}$ as
$$
\rY^{(\alpha)} = \rY^{(\alpha(0, p_{\text{split}}))} \cup \rY^{(\alpha(1,p_{\text{split}}))}
$$
and we update the set of indices
$$
A_k \leftarrow \left(A_k \setminus \{\alpha\}\right) \cup \{ \alpha(0, p_{\text{split}}), \alpha(1, p_{\text{split}}) \},
$$
i.e., we replace the splitted index $\alpha$ by the refined sets $\alpha(0, p_{\text{split}})$ and $\alpha(1, p_{\text{split}})$.

Once all the indices in $A_{k-1}$ have been visited (and splitted if needed), we finally set
\begin{align*}
\rY_k &= \cup_{\alpha \in A_k} \rY^{(\alpha)} \\
\rL_k &= \{\Valpha\}_{\alpha\in A_k}.
\end{align*}
If
$$
e_k \coloneqq \max_{\alpha\in A_{k}} \ealpha \leq \eps,
$$
the algorithm stops and we set $K=k$, otherwise we increase the index $k$ and repeat the inductive step.
\end{itemize}
\end{itemize}

This algorithm, referred to as $\algoY$ in what follows, produces a library $\cL\coloneqq\rL_K \in  \cL_{2^{K-1}, n_{\max}}$. That is, the library $\cL$ has cardinality $|\cL|\le 2^{K-1}$ and contains spaces of dimension at most $n_{\max}$. Moreover, the total number of snapshots generating the library is $\sum_{\alpha\in A_K}\dim(V^{(\alpha)})\le n_{\max}2^{K-1}$.

Once $\cL$ has been computed, we need to define the final approximation $\hat u_\cL$ introduced in equation \eqref{eq:library-approach}, namely
\begin{equation*}
\hat u_\cL : \rY \to \cL, \quad \hat u_\cL(y)\in V^{(\iota(y))}.
\end{equation*}
Thanks to the dyadic splittings in the parameter domain, the mapping $\iota$ has a very simple characterization, and it reads
$$
\iota(y) = \{\alpha \in A_k \cond y \in \rYalpha\}.
$$

\subsection{An approach which does not depend on the shape of $\rY$}
\label{sec:general-split}

The previous strategy heavily relies on the assumption that the parameter space $\rY$ has a tensor product structure. We introduce in this subsection an alternative algorithm which can be used to generate a library $\cL=\{\Valpha\}_{\alpha\in A}$ for a generic $\rY$. The main difference lies in the definition of the subdomain $\rYalpha\subseteq\rY$ and the space $V^{(\alpha)}$ associated to an index $\alpha\in A$. However, once defined, the subdomain $\rYalpha\subseteq\rY$ generates again the subset
$$
\cMalpha = \{ u(y)\cond y \in \rYalpha \}\subseteq \cM.
$$
Here, the indices $\alpha$ will be binary, namely of the form $\alpha=(a_1,\ldots,a_k)\in\{0,1\}^k$. The root of the tree is encoded with $\alpha=(1)$ and we work with the convention
$$
\cMalpha=\cM,\; \rYalpha=\rY,\quad \text{for }\alpha=(1).
$$
Moreover, for $u\in\cM$, we will use the notation
$$\dist(u,V^{(\alpha)})\coloneqq\inf_{v\in V^{(\alpha)}}\dist(u,v).$$

The second algorithm proposed to generate a tree-based library approximation can be formulated as follows. We set the target tolerance $\eps>0$. We then proceed iteratively as follows:
\begin{itemize}
\item Initialization ($k=1$): We define
\begin{equation*}
\alpha_k = (1), \qquad A_k = \{\alpha_k\}, \qquad \rY_k = \cup_{\alpha \in A_k} \rY^{(\alpha)} = \{ \rY \}.
\end{equation*}
When $V$ is a Hilbert or a Banach space, we select a first snapshot applying one step of the plain greedy algorithm, that is,
\begin{equation}
\label{eq:init-general-split}
[\rU^{(\alpha_k)}, V^{(\alpha_k)}, e^{(\alpha_k)}] = \greedy(\cM^{(\alpha_k)}, V, \eps, n_{\max}=1),
\end{equation}
where $\rU^{(\alpha_k)}=\{u_{(\alpha_k)}\}$ and where $\cM^{(\alpha_k)}=\cM$ by construction. We then set
$$
\rL_{k}=\{ V^{(\alpha_k)}\}.
$$
When $V=(\Pr,W_2)$, we apply the same strategy but we need to select two snapshots to initialize the greedy algorithm so we should use $n_{\max}=2$ in \eqref{eq:init-general-split}.

If $e^{(\alpha_k)}\leq \eps$, the algorithm stops, otherwise we set $k=2$ and the algorithm goes into the inductive step.
\item Inductive step $k>1$: From the previous step, we receive $A_{k-1},\,\rY_{k-1}$, and $\rL_{k-1}=\{\Valpha\}_{\alpha\in A_{k-1}}$. We set $A_{k}=A_{k-1}$, and we are going to update $A_k$ with the procedure that follows. For every $\alpha\in A_{k-1}$ for which $e^{(\alpha)}>\eps$, we split $\cMalpha$ into two subsets as follows. We make two steps of a restarted greedy algorithm that takes $\Valpha$ as a starting approximation space,
\begin{equation*}
[\widetilde \rU^{(\alpha)} , \widetilde \rV^{(\alpha)}, e]=
\restartedgreedy[\cM^{(\alpha)}, V, \Valpha, \eps, n_{\max}=\dim(\Valpha)+2].
\end{equation*}
By construction,
$$
\widetilde \rU^{(\alpha)} = \rUalpha \cup\{u_{(\alpha,0)},u_{(\alpha,1)}  \},
$$
where the two snapshots $u_{(\alpha,0)},u_{(\alpha,1)}$ are the ones that have been added by the restarted greedy algorithm. For $i\in \{0,1\}$, we define
\begin{equation*}
\rU^{(\alpha,i)} \coloneqq \rU^{(\alpha)}\cup \{u_{(\alpha,i)}\} \quad \text{and} \quad V^{(\alpha,i)}\coloneqq 
\begin{cases}
\vspan\{\rU^{(\alpha,i)}\}\quad\; \text{if}\; V=\text{Banach space} \\
\bary(\Sigma_r, \rU^{(\alpha,i)})\; \text{if}\; V=(\Pr,W_2)
\end{cases}
\end{equation*}
with $r\coloneqq|\rU^{(\alpha,i)}|$. Next we decompose $\rY^{(\alpha)} = \rY^{(\alpha,0)} \cup \rY^{(\alpha,1)} $ according to
$$
\rY^{(\alpha,0)} \coloneqq\left\{y\in \rY^{(\alpha)}: \quad \dist(u(y),V^{(\alpha,0)})\le \dist(u(y),V^{(\alpha,1)}) \right\}
$$
and $\rY^{(\alpha,1)}\coloneqq\rY^{(\alpha)}\setminus \rY^{(\alpha,0)}$. Then we define
$$e^{(\alpha,i)}\coloneqq \dist(\cM^{(\alpha,i)},V^{(\alpha,i)}), \quad i=0,1,$$
and we update
$$
A_k \leftarrow \left(A_k \setminus \{\alpha\}\right) \cup \{ (\alpha,0), (\alpha,1) \},
$$
where we have replaced the splitted index $\alpha$ by the refined sets $(\alpha,0)$ and $(\alpha,1)$. Once all the indices in $A_{k-1}$ have been visited (and splitted if needed), we finally set
\begin{align*}
\rY_k &= \cup_{\alpha \in A_k} \rY^{(\alpha)} \\
\rL_k &= \{\Valpha\}_{\alpha\in A_k}.
\end{align*}
If
$$
e_k \coloneqq \max_{\alpha\in A_{k}} \ealpha \leq \eps,
$$
the algorithm stops and we set $K=k$, otherwise we increase the index $k$ and repeat the inductive step.
\end{itemize}

This algorithm, referred to as $\algoM$ in what follows, produces a library $\cL=\rL_K \in  \cL_{2^{K-1}, n}$, where $n=K$ if $V$ is a Banach space and $n=K+1$ if $V=(\Pr,W_2)$. In this case, the total number of snapshots building the library corresponds to the number of nodes in the underlying tree, which may be much smaller than $\sum_{\alpha\in A_K}\dim(V^{(\alpha)})$. Indeed, many snapshots belong to several spaces, in particular the snapshot(s) associated to the root, i.e. $\rU^{(1)}$, belong to all the spaces in $\cL$.

We define the final approximation $\hat u_\cL$ as
\begin{equation*}
	\hat u_\cL : \rY \to \cL, \quad \hat u_\cL(y)\in V^{(\iota(y))},
\end{equation*}
where
\begin{equation}
\iota(y) = \{\alpha \in A_k \cond y \in \rYalpha\}
\label{eq:iota-algoM}
\end{equation}
is defined similarly to the previous splitting strategy involving cartesian splittings of the parameter domain. In contrast to that approach, here the mapping $\iota$ does not have a simple characterization in general because the shape of the resulting subdomains in $\rY$ becomes very general (see Figures \ref{fig:Diff_alpha_domains} and \ref{fig:diff_pw-conv} for an illustration). Among the possible strategies to build this mapping, we may cite classification techniques from machine learning. We can also consider strategies based on metric embeddings where we learn distances of elements from $\cM$ in $V$ through distances in $\rY$ (see \cite{DFM2023}).

\begin{remark}
Instead of making two steps in the restarted greedy algorithm, one can alternatively do more than two steps, and take the last two selected snapshots to perform the split into two subsets. In the numerical tests, this alternative is done for the KdV example considered in Section~\ref{sec:numerics:kdv} (the restarted greedy algorithm makes three steps). Moreover, other strategies than $\restartedgreedy$ could be used for the selection of the additional snapshots. For instance, when $V$ is a Hilbert space with inner product $(\cdot,\cdot)_V$, we could select all the elements for which the approximation error in $V^{(\alpha)}$ is a fraction (e.g., a half) of the largest error. Then among them, we pick the two that are the less correlated, namely $v$ and $w$ for which $(v,w)_V/(\|v\|_V\|w\|_V)$ is minimized. We have tested this approach for the diffusion problem of Section~\ref{sec:diff}, but the results are not reported since they are similar to the ones obtained with two steps of greedy.
\end{remark}

\begin{remark}
In the Hilbert setting, we can show that if the Kolmogorov $n$-width \eqref{eq:kolmo-width} of $\cM$ decays algebraically (resp. exponentially), then so does the error $\dist(\cM, \cL)$ of the library $\cL$ generated by $\algoM$. The proof follows \cite[Corollary 3.3]{DPW2013}, but extended to local spaces and accounting for the fact that snapshots in $\rU^{(\alpha)}$ may not belong to $\cM^{(\alpha)}$. In other words, we can show that $\algoM$ performs no worse than $\algoPG$, which is expected due to \eqref{rel:widths}. We leave further theoretical analysis of the performance of both $\algoM$ and $\algoY$ for future work.
\end{remark}

\subsection{Advantages and disadvantages of each method}
We discuss the advantages and limitations of each approach in the light of different aspects:
\begin{itemize}
\item \textbf{Stopping criteria:} A direct comparison of $\algoM$ and $\algoY$ is actually not trivial since their stopping criteria are not identical. Indeed, we only prescribe the target accuracy for $\algoM$ (without control on the dimension of the spaces unless we add a restriction on the depth of the tree), while for $\algoY$ we prescribe $\varepsilon$ and $n_{\max}$. As a result, we know in advance the maximum dimension of the local spaces and of course, the smaller $n_{\max}$ the more spaces in the library. One could consider modifications of $\algoY$ that would not involve $n_{\max}$. For example, in the spirit of $hp$-adaptivity, we could iterate with greedy algorithms over the splitted subdomains until the sum of their approximation errors becomes smaller than a fraction of the approximation error of the current parent. This option would be interesting to consider, but it has not been explored in our numerical examples.
\item \textbf{Offline stage:} The construction of a library with $\algoY$ is much more time consuming than with $\algoM$. This is mainly due to the fact that we are testing all possible splits and selecting the best one in the sense of \eqref{eq:p_split}.
\item \textbf{Online cost:} The numerical complexity to solve the parametric PDE for a given parameter $y\in \rY$ is related to the dimension $n(\alpha)$ of the local spaces $\Valpha$. For both library approaches, this dimension is expected to be significantly smaller that the one for the \algoPG~as soon as we deal with problems that are not elliptic. So even if we deal with a collection of spaces, both library approaches are expected to provide significant gains in the online phase if one needs to compute many queries.

In the Hilbert space setting, the online cost scales like $n(\alpha)^2$ (solution of a linear system with a dense matrix of size $n(\alpha)\times n(\alpha)$). For the Wasserstein space, evaluating Wasserstein barycenters with $n(\alpha)$ components scales proportionally to $n(\alpha)^2$ in dimension 1, and it becomes dependent on the number of discretization points $\cN$ for larger dimensions (order $n(\alpha)\cN^2\log\cN$ with state of the art algorithms based on entropic regularization \cite{PC2019}).
\item \textbf{The mapping $\iota$:} For a given $y\in\rY$, we need to decide which space of the library to use, namely construct a mapping $\iota:\rY\rightarrow A_k$. While this mapping is given by construction for \algoY, we need to build it for $\algoM$. This can be done by using the training snapshots with an additional approximation procedure in the offline stage. Among the several options that one could consider, we may mention strategies based on learning metric embeddings as was introduced in \cite{DFM2023}: we learn distances over elements in $\cM$ by learning Euclidean distances in the parameter space. Then, for a new parameter $y\in\rY$, we select the space to use by searching for the element of the training set which is at closest distance (and for which we know that it is approximated by a certain space from $\cL$). It would also be possible to train a classifier.
\end{itemize}

\section{Numerical Examples}
\label{sec:numerics}

We illustrate the behavior of the algorithms for PDEs of different nature. In Section \ref{sec:numerics:hilbert}, we consider manifolds $\cM$ connected to diffusion and convection-diffusion parametric PDEs defined on Hilbert spaces. In Section \ref{sec:numerics:kdv}, we consider a parametric Korteweg-De Vries equation, and we view the generated manifold $\cM$ as a subset of the Wasserstein space. We summarize the main results that emerge from these tests in Section \ref{sec:numerics:conclusions}, and we use the examples to discuss the advantages and disadvantages of each approach.

\subsection{Diffusion and convection-diffusion problems in $V=H^1_0(\Omega)$}
\label{sec:numerics:hilbert}

In this section, we set $\Omega=(0,1)$ and we consider the Hilbert space $V=H_0^1(\Omega)$ endowed with the norm  $\|v\|_V\coloneqq \|v'\|_{L^2(\Omega)}$, $v\in H_0^1(\Omega)$. All the problems can be written in the form: given $y\in\rY$, find $u(y)\in V$ such that
\begin{equation}
\label{eq:weak-formulation-num}
a(u(y),v;y) = F(v), \quad \forall\, v\in V,
\end{equation}
where $a(\cdot,\cdot;y)$ is a bilinear form parametrized by $y$ which is uniformly coercive and continuous with constants $0<r\le R<\infty$. We let $V_h$ be the space of piecewise linear ($\mathbb{P}_1$) finite elements based on a uniform partition of $\Omega$ with mesh size $h=2^{-12}$, and consider the Galerkin projection $u_h(y)=P_{V_h}u(y)\in V_h$ of $u(y)$ to be the \emph{truth solution}. Then, according to \eqref{eq:equiv_res_err}-\eqref{eq:weak_greedy_est} and given the current approximation space $V_n\subset V_h$, we use the residual
$$\Delta_n(y)\coloneqq \|\cR_n(\cdot;y)\|_{V_h'}, \quad \cR_n(v;y)\coloneqq  F(v)-a(P_{V_n}u(y),v;y)$$
instead of the error $\|u_h(y)-P_{V_n}u(y)\|_V$ to drive the greedy algorithms.  We use a (deterministic) training set $\rY^{(N)}\subset\rY$ of finite dimension $N$ for the greedy algorithms, and we \emph{randomly} select one point in $\rY^{(N)}$ to generate the first snapshot instead of \eqref{sec:vanilla-greedy}. Finally, we set the tolerance to $\varepsilon=10^{-6}$.

\subsubsection{Diffusion problems} \label{sec:diff}
We set $\rY=[-1,1]^2$ (case $p=2$) and consider the diffusion model problem
\begin{equation} \label{eqn:diff}
-\frac{\partial}{\partial x}\left(a(x,y)\frac{\partial u}{\partial x} (x,y)\right) = 1, \quad x\in(0,1), \quad u(0,y) = u(1,y) = 0.
\end{equation}
The problem has a weak formulation of the form \eqref{eq:weak-formulation-num} which reads: find $u\in H^1_0(\Omega)$ such that
\begin{equation*}
\int_\Omega a(x, y) \frac{\partial u}{\partial x} (x,y) \frac{\partial v}{\partial x}(x)\dx = \int_\Omega v(x)\dx, \quad \forall v\in H^1_0(\Omega).	
\end{equation*}
We next consider two examples corresponding to two different choices for the diffusion coefficients $a$.
\\

\begin{figure*}[htbp]
	\centering
	\begin{subfigure}[t]{0.5\textwidth}
		\centering
		\includegraphics[scale=0.35]{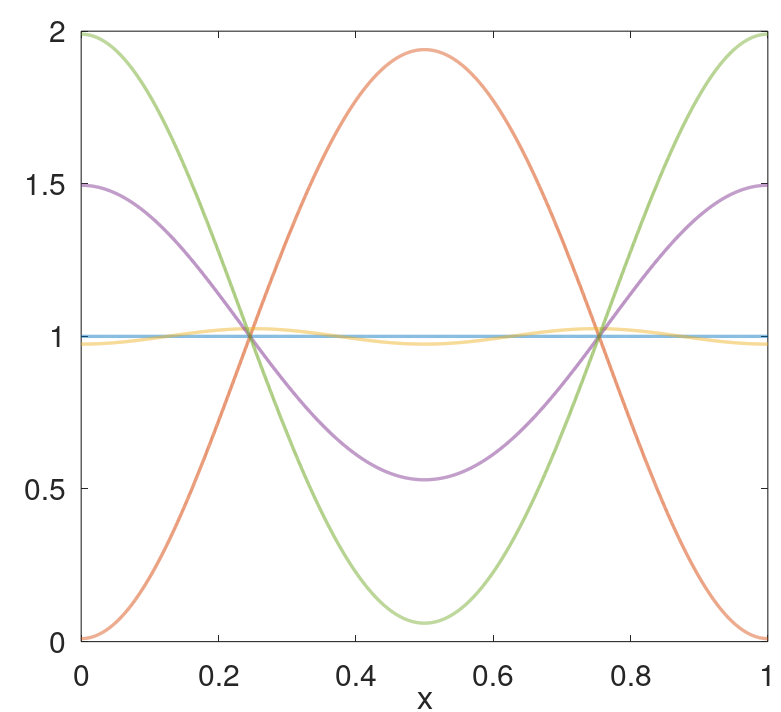}	
		\caption{Coefficient in \eqref{def:diff_coeff} with $\alpha=0.105$.}
		\label{fig:diff1_plot_a}
	\end{subfigure}%
	~ 
	\begin{subfigure}[t]{0.5\textwidth}
		\centering
		\includegraphics[scale=0.35]{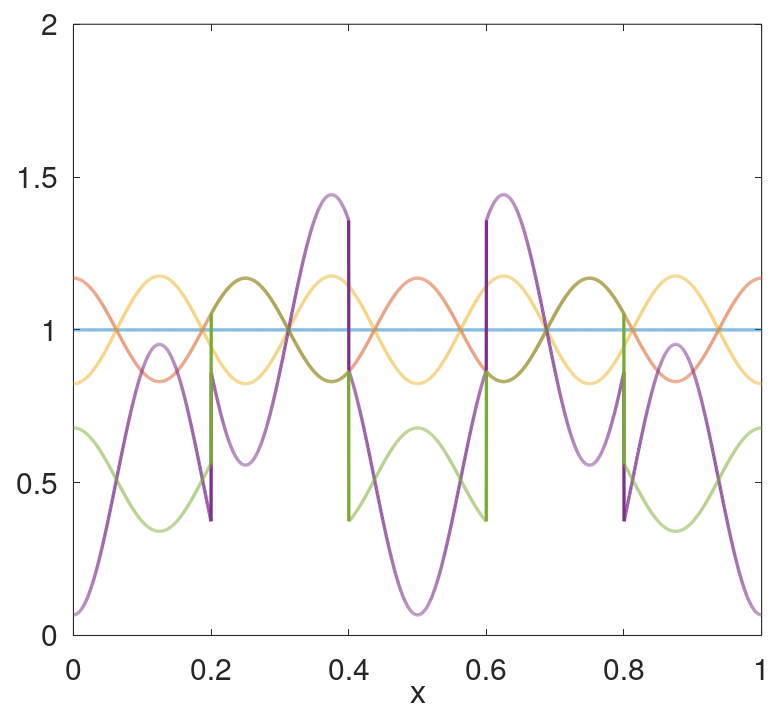}
		\caption{Coefficient characterized by \eqref{eqn:diff_a_pw_part1}-\eqref{eqn:diff_a_pw_part2}.}
		\label{fig:diff2_plot_a}
	\end{subfigure}
	\caption{\textbf{Diffusion problems:} Diffusion coefficient $a(x,y)$ for the following values of the parameter $y$: $(0,0)$, $(-1,-1)$, $(0,-1)$, $(0.5,0.5)$, $(1,1)$.}
\end{figure*}

\noindent \textbf{Test 1 (labeled as \diffone):} We first choose the diffusion coefficient $a$ defined for $y=(y_1,y_2)\in \rY$ and $x\in\Omega$ by
\begin{equation} \label{def:diff_coeff}
a(x,y)=1+\sum_{i=1}^2a_i(x)g_i(y_i) = 1+\frac{\cos(2\pi x)}{\alpha \pi^2}y_1+\frac{\cos(4\pi x)}{(2\pi)^2}y_2,
\end{equation}
where $\alpha>\alpha_{\min}=\frac{4}{4\pi^2-1}\approx 0.104$ so that the problem is coercive. We examine the impact of coercivity by examining the different manifolds $\cM$ generated by taking $\alpha=1$ and $\alpha=0.105$. The value $\alpha=0.105$ corresponds to a degeneration of coercivity in the sense that $r$, $R$ and the conditioning $\nu=R/r$ deviate from $1$ ($r=0.0097$, and $R=1.9903$), see Figure~\ref{fig:diff1_plot_a} for some realizations of $a(\cdot,y)$. The value $\alpha=1$ corresponds to a more benign coercive case with $r=0.8733$ and $R=1.1267$. The training set is $\rY^{(N)}$ with $N=40000$ corresponding to the tensor product grid with $200$ equally spaced points in each direction.

\begin{figure}[htbp]
	\centering
	\includegraphics[scale=0.35]{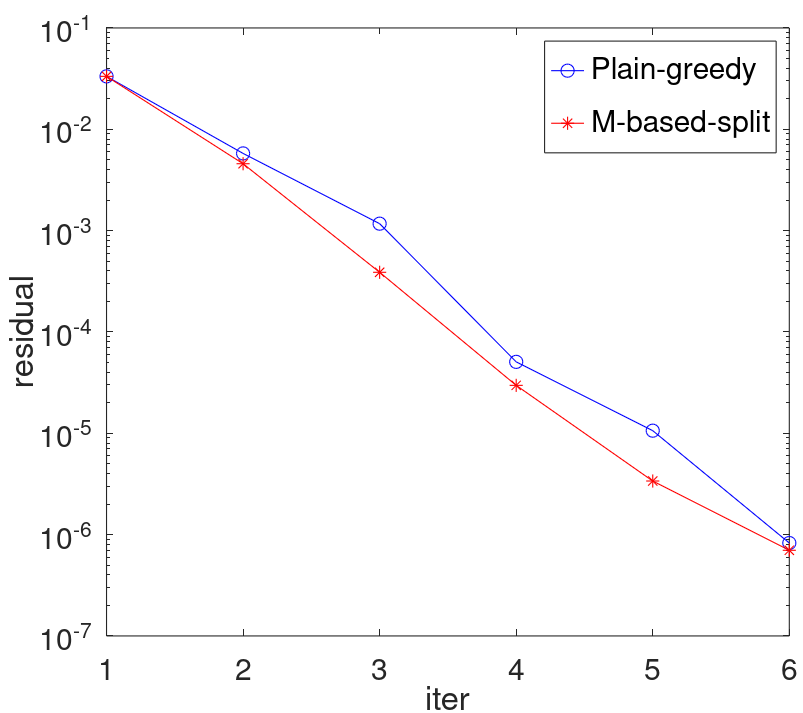}
	\hspace*{0.5cm}\includegraphics[scale=0.35]{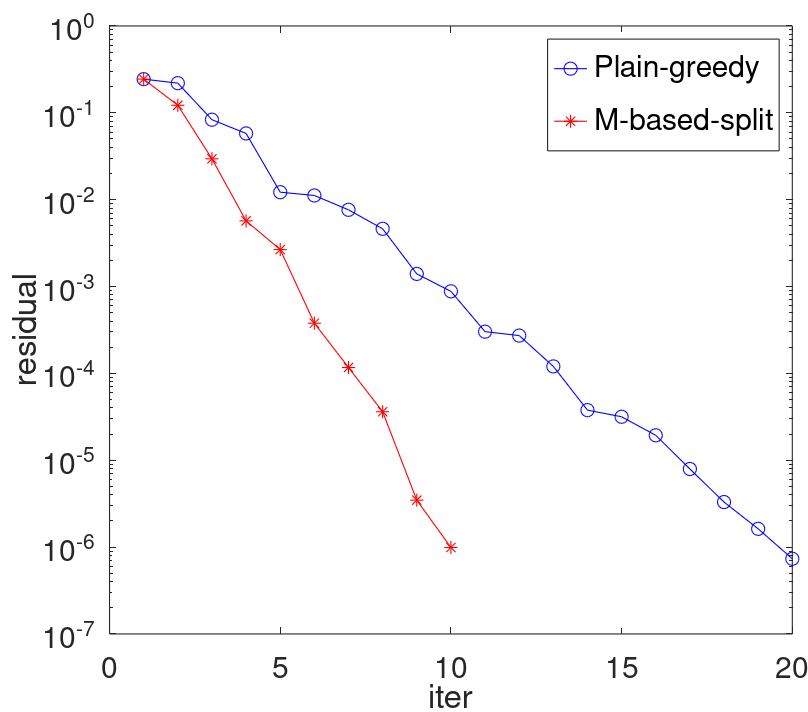}	
	\caption{\textbf{\diffone~problem:} Comparison of the algorithms $\algoPG$ and $\algoM$ with $\alpha=1$ (left) and $\alpha=0.105$ (right).}
	\label{fig:diff_alpha-conv}
\end{figure}

\begin{figure}[htbp]
	\centering
	\includegraphics[scale=0.5]{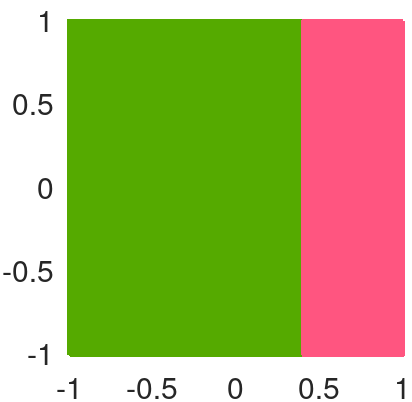}
	\hspace*{0.2cm}
	\includegraphics[scale=0.5]{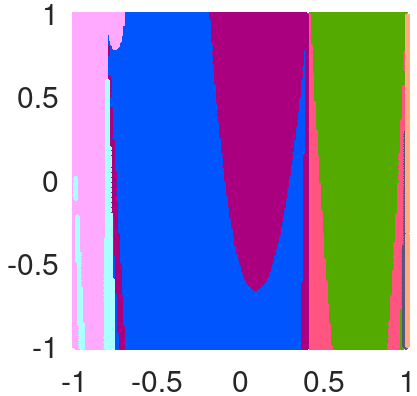}		
	\hspace*{0.2cm}
	\includegraphics[scale=0.5]{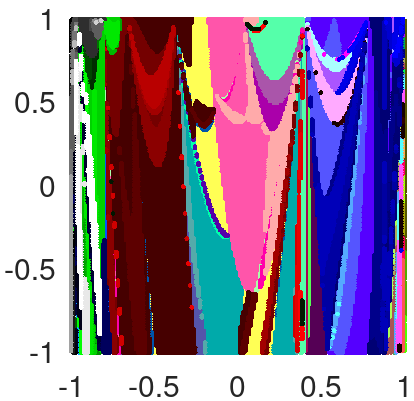}
	\caption{\textbf{\diffone~problem:} Partition of $\rY^{(N)}$ obtained with algorithm $\algoM$ after 2 iterations (left), 4 iterations (middle) and 10 iterations (right) in the case $\alpha=0.105$.}
	\label{fig:Diff_alpha_domains}
\end{figure}

When $\alpha=1$, the algorithm $\algoPG$ converges in $6$ iterations while the algorithm $\algoM$ generates a library with $11$ spaces of dimension at most $6$ ($5$ spaces of dimension $5$ and $6$ of dimension $6$), see Figure~\ref{fig:diff_alpha-conv}-left for the evolution of the maximum residual as a function of the iteration $k$ (at iteration $k$, all the RB spaces are of dimension at most $k$). Therefore, when coercivity is not degenerate there is no advantage of using the nonlinear reduced model. The situation is different when $\alpha \to \alpha_{\min}$ and coercivity becomes degenerate. We illustrate this by considering the manifold $\cM$ generated when $\alpha=0.105$. In this case, we need a space of dimension $n=20$ to reach $\eps$ with the \algoPG~method while the maximum dimension of the space in the library is $10$, see Figure~\ref{fig:diff_alpha-conv}-right. More precisely, the library contains $96$ spaces. It is composed of $\{3,6,31,37,19\}$ spaces of dimension $\{6,7,8,9,10\}$. The spaces of the library are generated using a total of $261$ snapshots. The partition of the parameter domain obtained at different steps is given in Figure~\ref{fig:Diff_alpha_domains}.

\begin{figure}[htbp]
	\centering
	
	\begin{tikzpicture}[scale=1.5]
	
	\filldraw[fill=blue!30, thick] (0,-1) -- (0.5,-1) -- (0.5,1) -- (0,1) -- (0,-1);
	
	\filldraw[fill=green!30, thick] (-0.5,-1) -- (0,-1) -- (0,0) -- (-0.5,0) -- (-0.5,-1);

	\filldraw[fill=red!30, thick] (-1,0.5) -- (0,0.5) -- (0,1) -- (-1,1) -- (-1,0.5);

	\filldraw[fill=red!30, thick] (0.5,0) -- (1,0) -- (1,1) -- (0.5,1) -- (0.5,0);

	\filldraw[fill=green!30, thick] (-0.75,-1) -- (-0.5,-1) -- (-0.5,0) -- (-0.75,0) -- (-0.75,-1);

	\filldraw[fill=red!30, thick] (-1,0) -- (-0.5,0) -- (-0.5,0.5) -- (-1,0.5) -- (-1,0);	
	
	\filldraw[fill=orange!30, thick] (-0.5,0) -- (0,0) -- (0,0.5) -- (-0.5,0.5) -- (-0.5,0);	
	
	\filldraw[fill=green!30, thick] (0.5,-1) -- (0.75,-1) -- (0.75,0) -- (0.5,0) -- (0.5,-1);
		
	\filldraw[fill=red!30, thick] (0.75,-1) -- (1,-1) -- (1,0) -- (0.75,0) -- (0.75,-1);
	
	\filldraw[fill=green!30, thick] (-0.875,-1) -- (-0.75,-1) -- (-0.75,0) -- (-0.875,0) -- (-0.875,-1);
	
	\filldraw[fill=red!30, thick] (-1,-1) -- (-0.93755,-1) -- (-0.9375,0) -- (-1,0) -- (-1,-1);			

	\filldraw[fill=green!30, thick] (-0.9375,-1) -- (-0.875,-1) -- (-0.875,0) -- (-0.9375,0) -- (-0.9375,-1);	
	
	\draw (-1,-0.95) node[left] {$-1$};
	\draw (-0.95,-1.05) node[below] {$-1$};
	\draw (1,-1.05) node[below] {$1$};
	\draw (-1,0.95) node[left] {$1$};
				
\end{tikzpicture}	
	
\caption{\textbf{\diffone~problem:} Final partition of $\rY^{(N)}$ generated by the algorithm $\algoY$ with $n_{\max}=10$ when $\alpha=0.105$. The local spaces in the library are of dimension $6$ (orange), $7$ (green), $8$ (blue) and $10$ (red).}
\label{fig:Diff_alpha_Cartesian}
\end{figure}
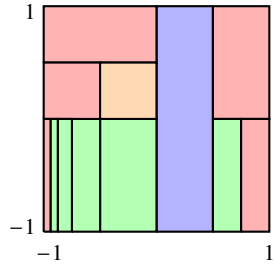

For a comparison, we run the algorithm $\algoY$ with $n_{\max}=10$. This algorithm produces a library with 12 spaces, see Figure~\ref{fig:Diff_alpha_Cartesian} for the partition of the domain and the dimension of the local RB spaces. A summary of the main features of the reduced models for this problem is provided in Table~\ref{numerics:maintable} (see label \diffone).

\begin{figure}[htbp]
	\centering
	\includegraphics[scale=0.35]{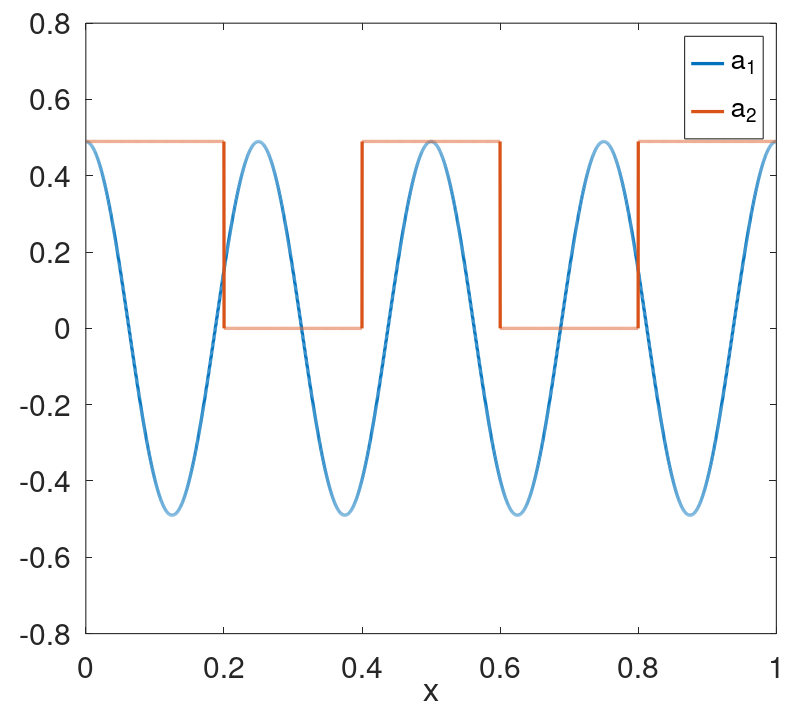}
	\hspace*{0.5cm}\includegraphics[scale=0.35]{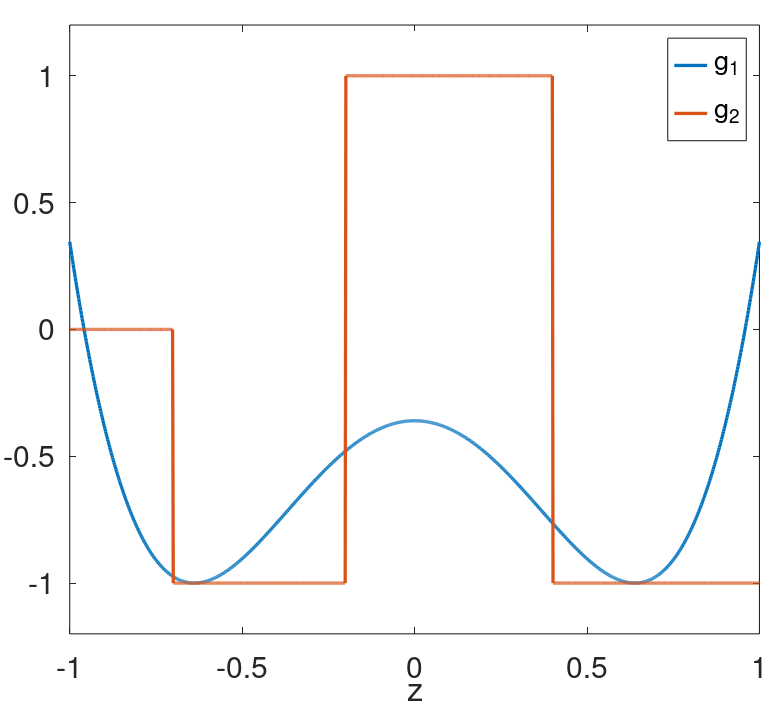}	
	\caption{\textbf{\difftwo~problem:} Plot of the functions defined \eqref{eqn:diff_a_pw_part1} and \eqref{eqn:diff_a_pw_part2}.}
	\label{fig:diff2_fcts}
\end{figure}

\noindent \textbf{Test 2 (labeled as \difftwo):}
We next consider a different diffusion coefficient, namely $a$ as in \eqref{def:diff_coeff} but with
\begin{equation} \label{eqn:diff_a_pw_part1}
a_1(x) = 0.49\cos(8\pi x), \quad g_1(y_1)=((1.4y_1)^2-0.8)^2-1	
\end{equation}
and
\begin{equation} \label{eqn:diff_a_pw_part2}
a_2(x) = 0.49\chi_{A_1}(x), \quad g_2(y_2)=\chi_{B_1}(y_2)-\chi_{B_2}(y_2),
\end{equation}
where $\chi$ denotes the indicator function and where $A_1=[0,0.2)\cup(0.4,0.6)\cup(0.8,1]$, $B_1=(-0.2,0.4)$ and $B_2=(-0.7,-0.2)\cup(0.4,1]$. We refer to Figure~\ref{fig:diff2_fcts} for a plot of the functions $a_1(x)$, $a_2(x)$, $g_1(z)$ and $g_2(z)$ for $x\in[0,1]$ and $z\in[-1,1]$.

\begin{figure}[htbp]
	\centering	
	\begin{subfigure}{.5\textwidth}
		\centering
		\includegraphics[scale=0.35]{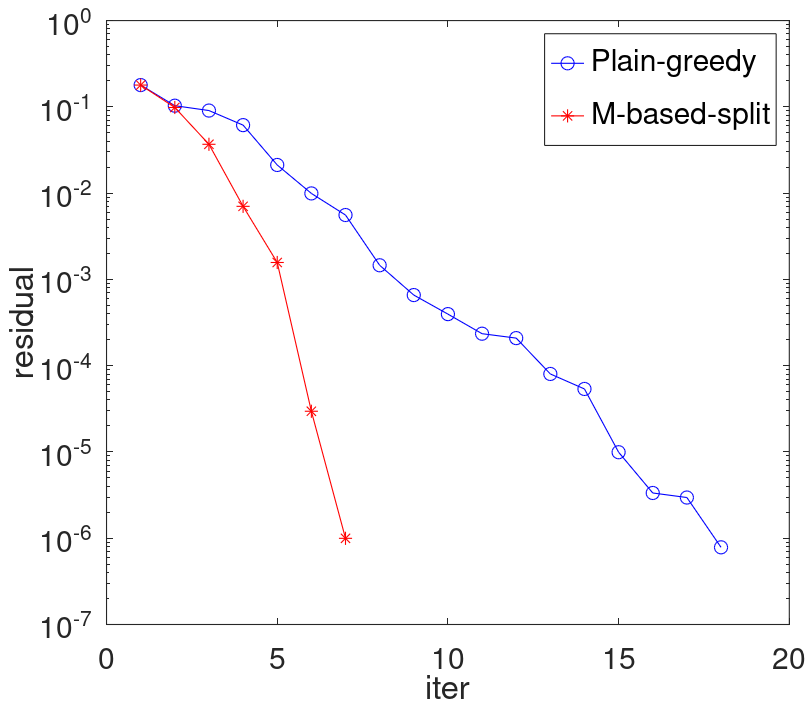}
	\end{subfigure}%
	\begin{subfigure}{.5\textwidth}
		\centering
		\includegraphics[scale=0.55]{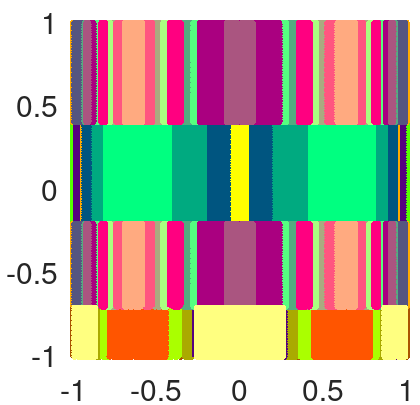}
	\end{subfigure}
	\caption{\textbf{\difftwo~problem:} Comparison of the algorithms $\algoPG$ and $\algoM$ (left) and final partition of $\rY^{(N)}$ generated by $\algoM$ (right).}
	\label{fig:diff_pw-conv}
\end{figure}

In this case, the algorithm $\algoPG$ generates a space of dimension $18$, while the library obtained with $\algoM$ has 24 spaces of dimension at most 7 ($\{3,11,10\}$ spaces of dimension $\{5,6,7\}$). The spaces of the library are generated with 61 snapshots. The maximum residual at each step of the algorithm is given in Figure~\ref{fig:diff_pw-conv}, which also contains the final partition associated to the library.

\begin{figure}[htbp]
	\centering	
	\includegraphics[scale=0.5]{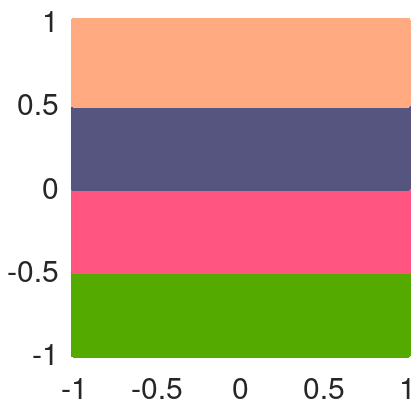}
	\hspace*{0.2cm}
	\includegraphics[scale=0.5]{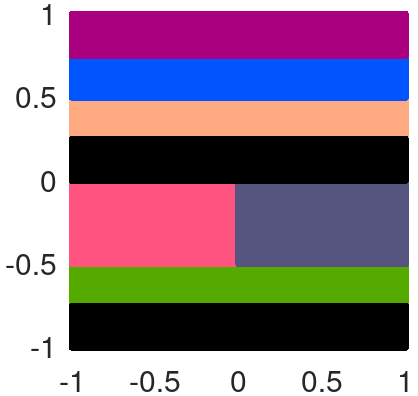}
	\hspace*{0.2cm}	
	\includegraphics[scale=0.5]{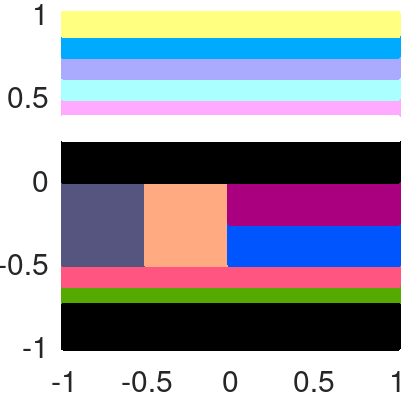}			
	\caption{\textbf{\difftwo~problem:} Partition of $\rY^{(N)}$ obtained with algorithm $\algoY$ after 3 iterations (left), 4 iterations (middle) and 5 iterations (right) when $n_{\max}=7$. The target accuracy is reached on the subdomains in black and in white.}
	\label{fig:diff2_Ysplit}
\end{figure}

For this problem, the agorithm $\algoY$ with $n_{\max}=7$ does not converge when we use \eqref{eq:p_split} to determine which direction to split at each step. The direction $y_2$ is chosen most of the times, see Figure~\ref{fig:diff2_Ysplit} for the partition of the parameter domain after a few steps of the algorithm. For instance, for $\rY^{(\alpha)}=[-1,1]\times[0.75,1]$, a split in $y_1$ gives the error $e_1^{(\alpha)}=\max\{2.2782\cdot 10^{-4},9.6067\cdot 10^{-5}\}$ while a split in $y_2$ gives the error $e_2^{(\alpha)}=\max\{2.2605\cdot 10^{-4},2.2605\cdot 10^{-4}\}$, and thus $p_{\text{split}}=2$. However, the error on the two subdomains thereby generated is not decreased. If we set instead $n_{\max}=10$ or we alternate between vertical ($p_\text{split}=1$) and horizontal ($p_\text{split}=2$) splits (starting with $y_1$, then $y_2$, $y_1$ and so on), then we obtain the results provided in Table~\ref{numerics:maintable} (see label \difftwo).

\begin{table}[htbp]
\begin{center}
	
\begin{tabular}{ c c c c c c c c}
& Algorithm & $|\cL|$ & Dim.~spaces       & Nb.~spaces        & Nb.~snap. & $\eps$  &$n_{\max} $\\ [2ex]
\hline
\multirow{3}{1cm}{\textbf{\diffone} ($\alpha=0.105$)} &
$\algoPG$   & 1       & \{20\}                     & \{1\}                   & 20          & $10^{-6}$ & $\infty$  \\  [2ex]
		
& \algoM    & 96$\;$  & $\{6,7,8,9,10 \}$  & $\{3,6,31,37,19\}$ & 261        & $10^{-6}$ & $\infty$   \\[2ex]

& $\algoY$    & 12$\:$      &  $\{6,7,8,10\}$     &  $\{1,5,1,5\}$   &  99  & $10^{-6}$ & $10$\\
\hline    	                 
\multicolumn{1}{l}{\textbf{\difftwo}} &		
$\algoPG$   & 1     & \{18\}                     & \{1\}                   & 18      & $10^{-6}$ & $\infty$      \\  [2ex]
 & $\algoM$    & 24$\;$  & $\{5,6,7 \}$  & $\{3,11,10\}$ & 61      & $10^{-6}$ & $\infty$     \\[2ex]
			
& $\algoY$    & 48 &  $\{8,9,10\}$     &  $\{6,36,6\}$   &  432  & $10^{-6}$ & $10$ \\[1ex]
& $\algoY^*$    & 314$\:$      &  $\{3,4,5,6,7\}$     &  $\{22,92,50,28,122\}$   &  1706  & $10^{-6}$ & $7$ \\  
\hline
\multicolumn{1}{l}{\textbf{\convdiff}} & $\algoPG$   & 1     & \{27\}                     & \{1\}                   & 27    & $10^{-6}$ & $\infty$        \\  [2ex]

& \algoM    & 42$\;$  & $\{5,6,7,8 \}$  & $\{3,11,10\}$ & 61    & $10^{-6}$ & $\infty$       \\[2ex]

& $\algoY$    & 11$\:$      &  $\{6,7\}$     &  $\{2,9\}$   &  75 & $10^{-6}$ & $8$  \\

\hline
\multicolumn{1}{l}{\textbf{\kdv}}
& $\algoPG$   & 1       & \{34\}                     & \{1\}                   & 34    & $10^{-3}$ & $\infty$        \\  [2ex]

 & $\algoM$    & 59$\;$  & $\{10, 12, 14, 16 \}$  & $\{2, 12, 19, 26\}$ & 208      & $10^{-3}$ & $\infty$     \\[2ex]

 & $\algoY$    & 20$\:$      &  $\{10,13,14,15\}$     &  $\{1,2,10,7\}$   &  281  & $10^{-3}$ & $15$\\
\hline  	 		
\end{tabular}
\caption{\textbf{Main features of the libraries obtained with the three different algorithms and for all the numerical examples}. Here, $|\cL|$ is the cardinality of the library, `Dim.~spaces' gives the dimensions of the spaces contained in $\cL$ and `Nb.~spaces' gives how many spaces of each dimension there are, `Nb.~snap.' is the amount of snapshots required to build all the spaces within $\cL$, $\eps$ is the target accuracy of the algorithms, and $n_{\max}$ is the prescribed maximal dimension of the spaces. The results for $\algoY^*$ correspond to the case where we alternate between vertical and horizontal splits.}
\label{numerics:maintable}
\end{center}
\end{table}

\vspace*{0.2cm}
\subsubsection{Convection-diffusion problem}

\begin{figure}[htbp]
	\centering
	\includegraphics[scale=0.35]{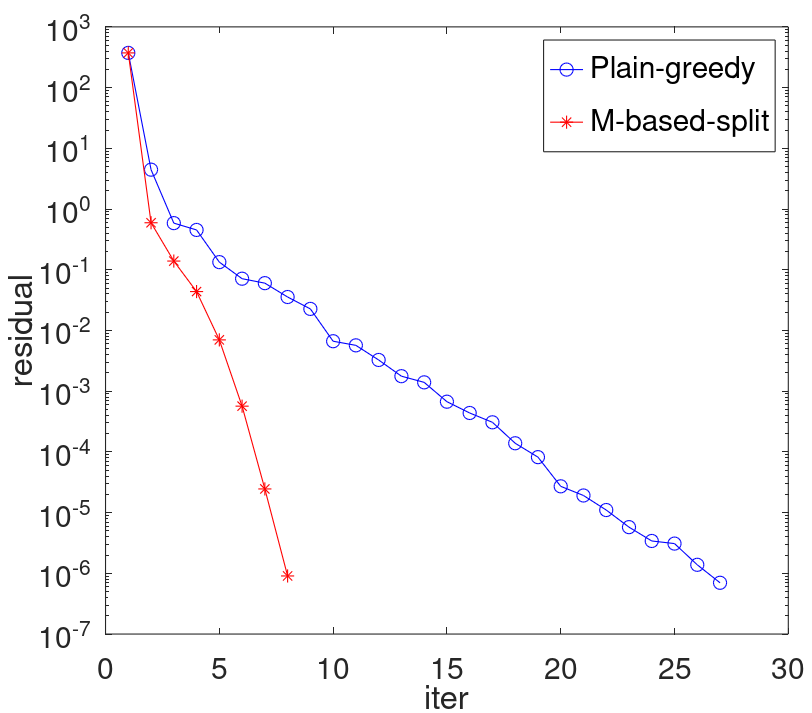}	
	\includegraphics[scale=0.35]{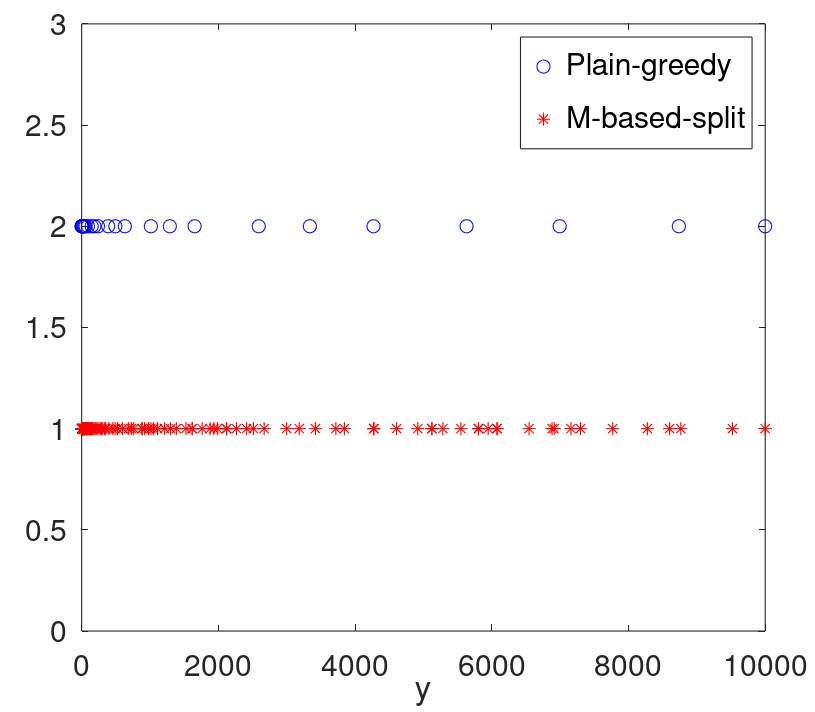}
	\caption{\textbf{\convdiff problem:} Comparison of the algorithms $\algoPG$ and $\algoM$.}
	\label{fig:diff_conv}
\end{figure}

Here $\rY=[y_{\min},y_{\max}]=[0,10000]$ and $u=u(x,y)$ is the (weak) solution to
\begin{equation} \label{eqn:conv_diff}
-\frac{\partial^2u}{\partial x^2}+y\frac{\partial u}{\partial x}=1 \quad \text{in } \Omega, \quad u(0,y)=u(1,y)=0.
\end{equation}
For this problem we have $r=1$ and $R=1+C_Py_{\max}$, where $C_P$ denotes the Poincar\'e constant in $\|v\|_{L^2(\Omega)}\le C_P\|v'\|_{L^2(\Omega)}$. Therefore, the condition number is such that $\nu=R/r\gg 1$ and we expect that the plain RB method will perform poorly (we refer to \cite{DPW2014} for a strategy developed to keep the condition number close to $1$). Here the training set $\rY^{(N)}$, $N=10001$, consists of the integers from $y_{\min}$ to $y_{\max}$. We give in Figure~\ref{fig:diff_conv} the evolution of the maximum residual at each step of the algorithm together with the location of all the selected parameters used to generate the final reduced model.

The algorithm $\algoPG$ converges in $27$ iterations, while the algorithm $\algoM$ generates a library with $42$ spaces of dimension at most $8$ ($4$ spaces of dimension $5$, $4$ of dimension $6$, and $5$ of dimension $7$) with a total of $111$ snapshots. If we run the algorithm $\algoY$ with $n_{\max}=8$, then the algorithm produces a library with 11 spaces, 2 of dimension 6 and 9 of dimension 7 and thus a total of 75 snapshots. Moreover, following the notation \eqref{eq:indices_alpha}, the partition of the parameter domain is $\rY=\cup_{\alpha\in A}\rY^{(\alpha)}$ with
$$\rY^{(\alpha)}\coloneqq y_{\max}\left[\frac{i}{2^k},\frac{i+1}{2^k}\right], \quad \alpha=(k,i), \qquad A = \{(k,1)\}_{k=1}^{10}\cup\{(10,0)\}.$$
If we set instead $n_{\max}=5$, then we obtain 22 spaces, all of dimension $5$, with the partition corresponding to $A = \{(k,2),(k,3)\}_{k=2}^{10}\cup\{(11,i)\}_{i=0}^3$. The results obtained for this problem are summarized in Table~\ref{numerics:maintable} (see label \convdiff).

\subsection{A KdV equation in $V=(\cP_2(\Omega), W_2)$ }
\label{sec:numerics:kdv}
We consider the same example from \cite{ELMV2020}, which is a two-soliton solution of the Korteweg-de-Vries equation (KdV). Expressed in normalized units, the PDE reads for all $x\in\bR$,
\begin{align*}
\partial_t u + 6u \partial_x u + \partial_x^3 u = 0.
\end{align*}
The equation admits a general $2$--soliton solution
\begin{eqnarray}
u(x,t) = -2 \partial^2_x \log \det(I + A(x,t)),
\end{eqnarray}
where $A(x,t) \in \bR^{2\times 2}$ is the interaction matrix whose components $a_{i,j}$ are
$$
a_{i,j}(x,t)= \frac{c_i c_j}{k_i + k_j} \exp\left( (k_i+k_j)x - (k_i^3 + k_j^3)t\right),\quad 1\leq i,j\leq 2.
$$
For any $t>0$, the total mass is equal to
$$
\int_\bR u(t,x) \dx = 4 (k_1+k_2).
$$
To illustrate the performance of our approaches, we set
$$
T=2.5\cdot 10^{-3},\;
k_1 = 30 - k_2.
$$
The parameter domain is\footnote{In \cite{ELMV2020}, $c_1$ and $c_2$ are fixed and only two parameters $(t, k_2)$ change. Here, we work with $c_1$ and $c_2$ as additional parameters.}
$$
\rY = \{ (t, c_1, c_2, k_2)\in [0, 2.5.10^{-3}] \times [0.9, 1.1] \times [0.2, 0.4] \times  [16, 22] \}
$$
It follows that the total mass is equal to $m=120$ for all the parametric solutions. We thus consider the solution manifold
$$
\cM = \{ \nu(y)(x) = \frac{u(y)}{m} \text{Leb(\dx)} \cond y\in \rY  \} \subset \Pr,
$$
where we view the solutions as probability density functions dominated by the Lebesgue measure (denoted as $\text{Leb(\dx)}$).

\begin{figure}[htbp]
\centering
\includegraphics[scale=0.4]{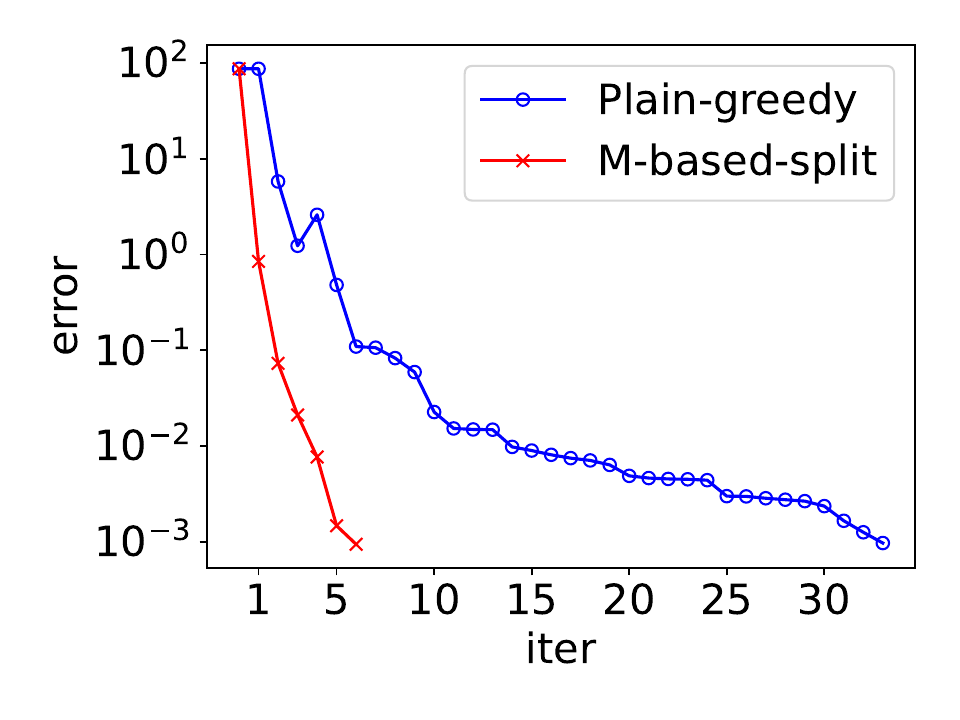}
\caption{\textbf{KdV problem:} Comparison of the algorithms $\algoPG$ and $\algoM$.}
\label{fig:kdv-conv-plain-greedy}
\end{figure}

We set $\eps=10^{-3}$ as the final target accuracy, and run the 3 algorithms with the same discrete set $\cM^{(N)}$ composed of $N=2\cdot 10^3$ training snapshots. The \algoPG~barycentric algorithm from Section \eqref{sec:barycentric-greedy} converges in $K=34$ iterations, and Figure \ref{fig:kdv-conv-plain-greedy} records the convergence of $\dist(\cM, V_n)$ as $n$ increases. Here, $V_n = \bary(\Sigma_n, \rU_n)$ is the barycenter space obtained with that algorithm.

\begin{figure}
\centering
\includegraphics*[scale=0.3]{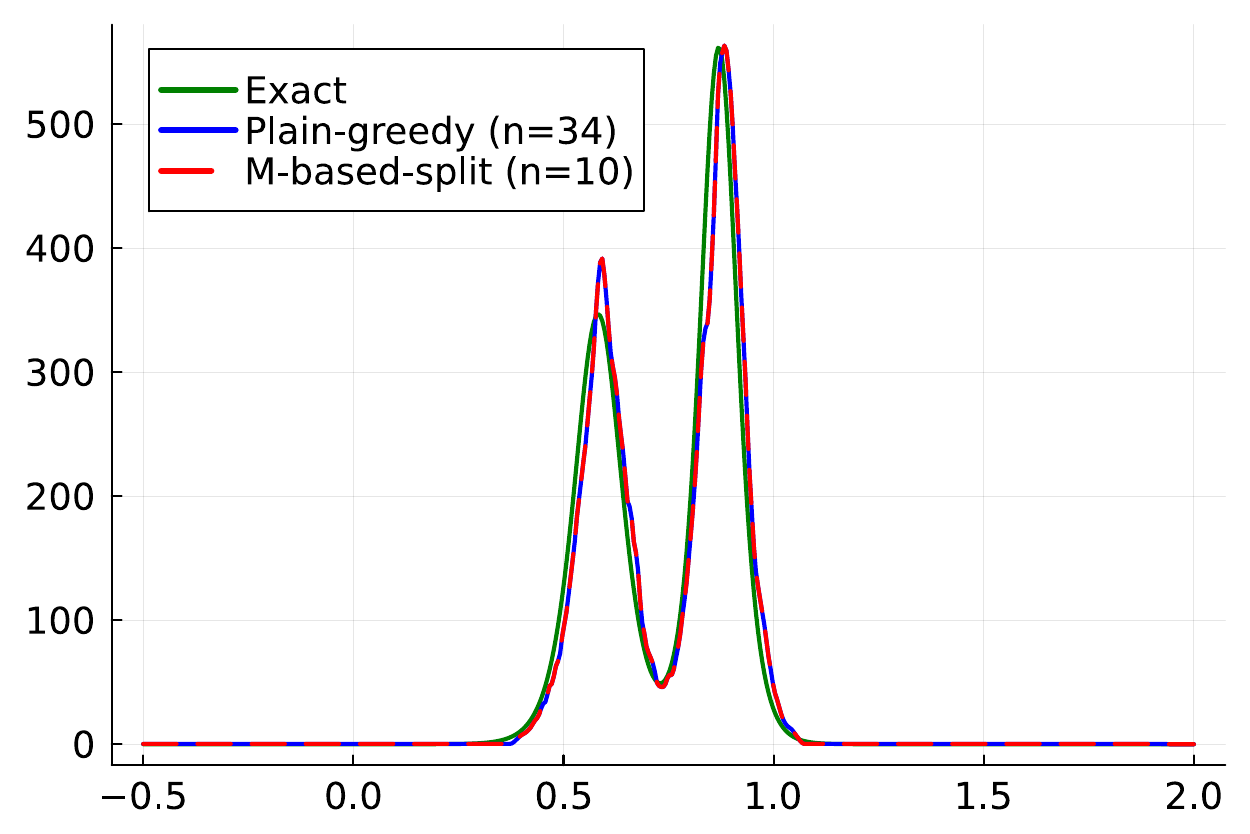}
\caption{\textbf{KdV problem:} Approximation of $u(y)$ for
$y=(2.23.10^{-3},  0.901, 0.386, 13.1)$ with $\algoPG$ and $\algoM$ once the stopping criterion as been reached. Visually, both approximations look identical and very close to the exact solution. However, \algoM~works with a space of smaller dimension ($n=10$ vs.~$n=34$).}
\label{fig:kdv-visualization}
\end{figure}

\begin{figure}[htbp]
\centering
\includegraphics[scale=0.2]{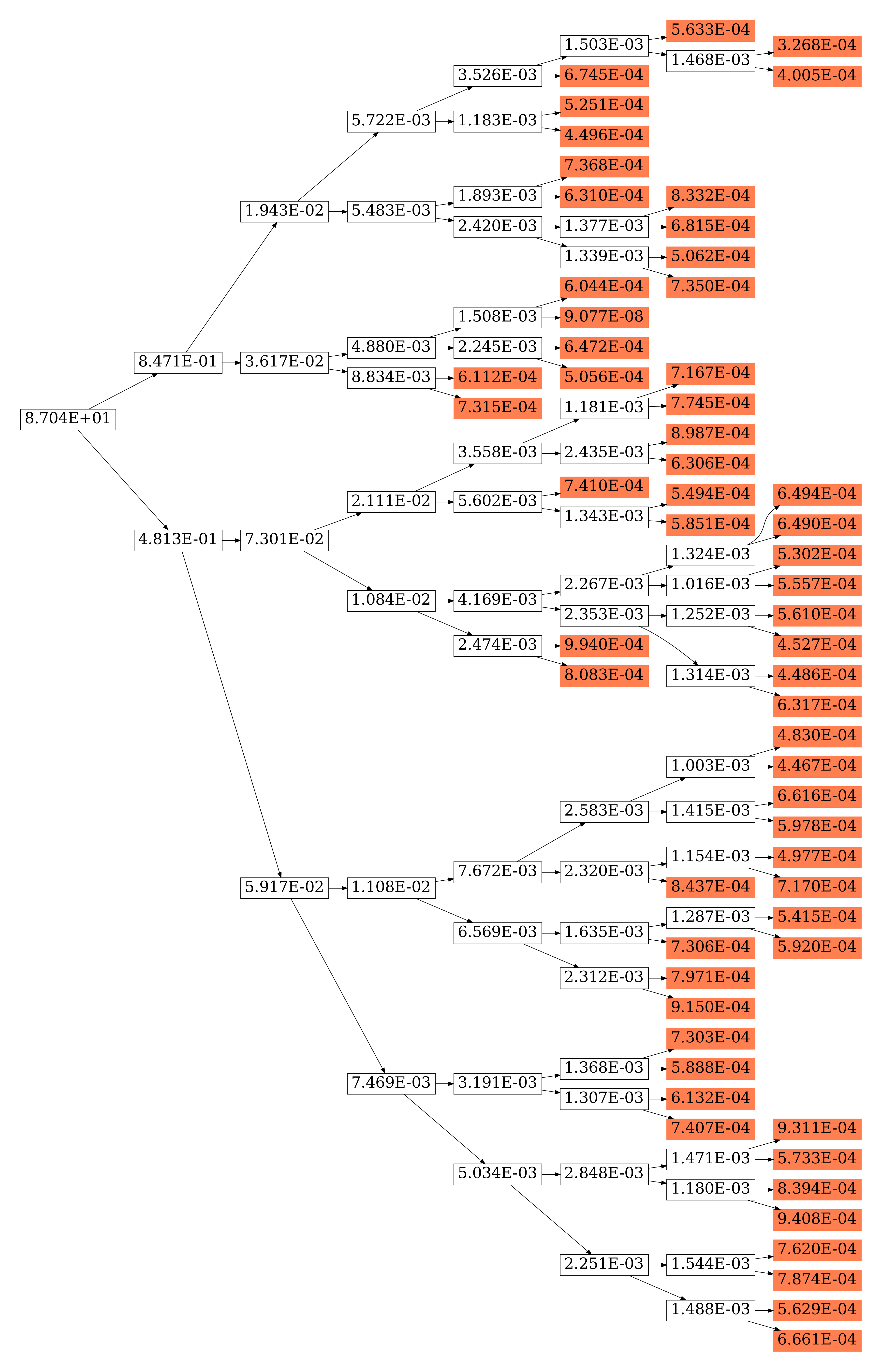}
\caption{\textbf{KdV problem:} Algorithm \algoM: accuracy $\dist(\cMalpha, \Valpha)$ reached at each node. The spaces $\Valpha$ that generate the library are located at the leaves of the tree (which  are colored in orange).}
\label{kdv:algoM-accuracy}
\end{figure}

\begin{figure*}[htbp]
\centering
\begin{subfigure}[t]{0.45\textwidth}
\centering
\includegraphics[width=\textwidth]{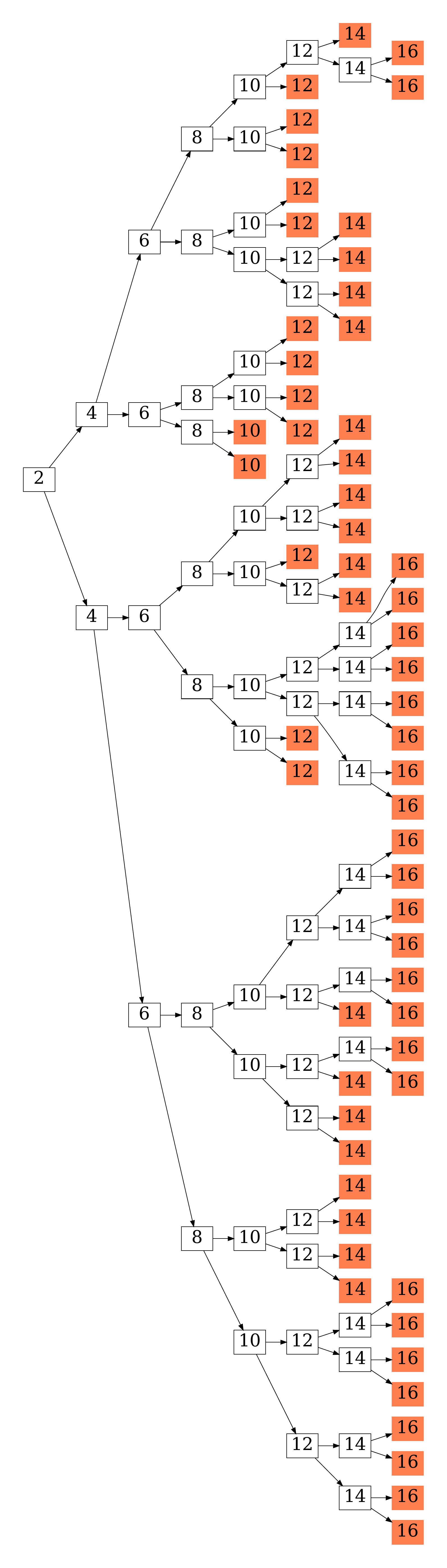}
\caption{Dimension $n(\alpha)=\dim(\Valpha)$ for each node.}
\label{kdv:algoM-dim}
\end{subfigure}%
~ 
\begin{subfigure}[t]{0.55\textwidth}
\centering
\includegraphics[width=\textwidth]{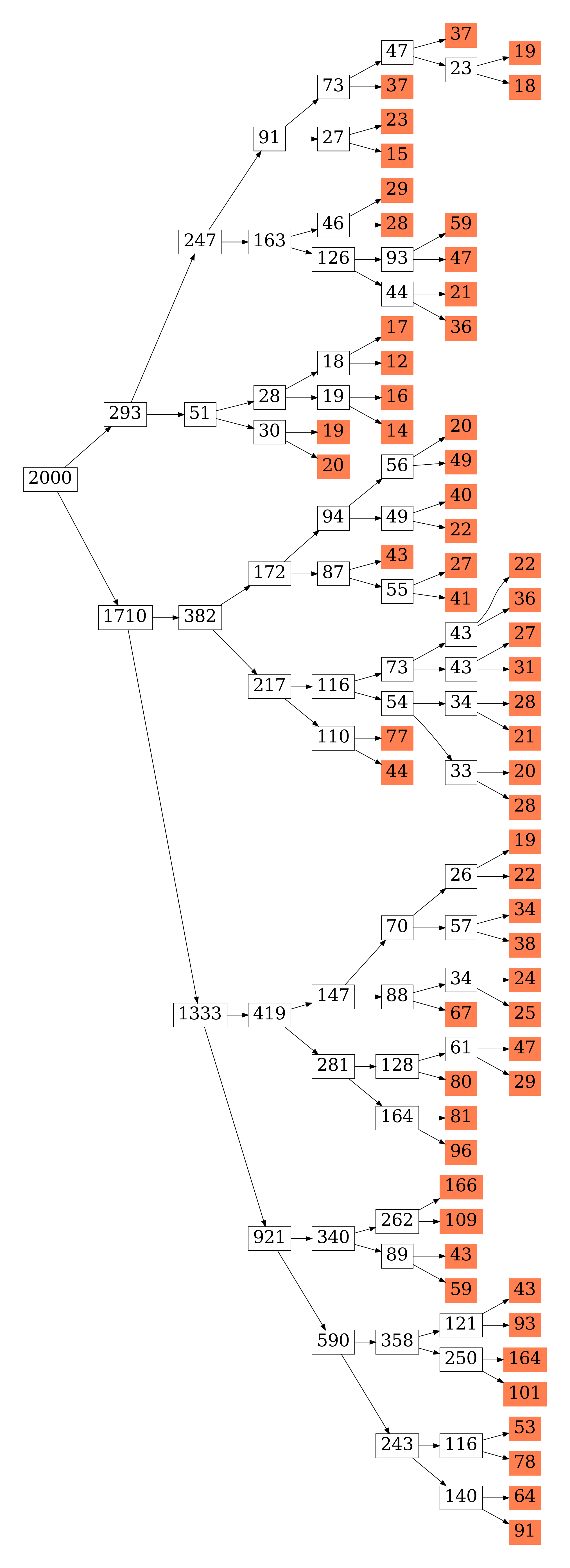}
\caption{Number of training snapshots for each node.}
\label{kdv:algoM-nb-snapshots}
\end{subfigure}
\caption{\textbf{KdV problem:} Algorithm $\algoM$: dimension and number of training snapshots in each tree node (leaves are colored in orange).}
\end{figure*}

We next comment on the behavior of the \algoM~algorithm. It converges in 6 iterations when we makes three steps in the restarted greedy algorithm. Figure \ref{fig:kdv-conv-plain-greedy} summarizes the approximation errors attained at each level of the tree. Figure \ref{kdv:algoM-accuracy} gives a similar type of information but in expanded form: it shows the intermediate accuracies that the algorithm reaches at each node. Figures \ref{kdv:algoM-dim} and \ref{kdv:algoM-nb-snapshots} show the space dimension $n(\alpha)=\dim(\Valpha)$, and the number of training snapshots at each node. The output of this algorithm is a library with $59$ spaces whose dimensions are $\{10, 12, 14, 16 \}$ and the amount of spaces for each dimension is $\{2, 12, 19, 26\}$. The spaces are the ones that are located at the leaves of tree. The $59$ spaces of the library are generated with $208$ snapshots, and some of them are shared by several different spaces. The main features of the library are summarized in Table \ref{numerics:maintable} (see label \kdv). Figure \ref{fig:kdv-visualization} shows the approximation of a randomly selected snapshot (that was not used in the training set) with \algoPG~and \algoM. Visually, both approximations look identical and very close to the exact solution. However, \algoM~works with a space of smaller dimension ($n=10$ vs.~$n=34$).

\begin{figure}[htbp]
\centering
\includegraphics[scale=0.2]{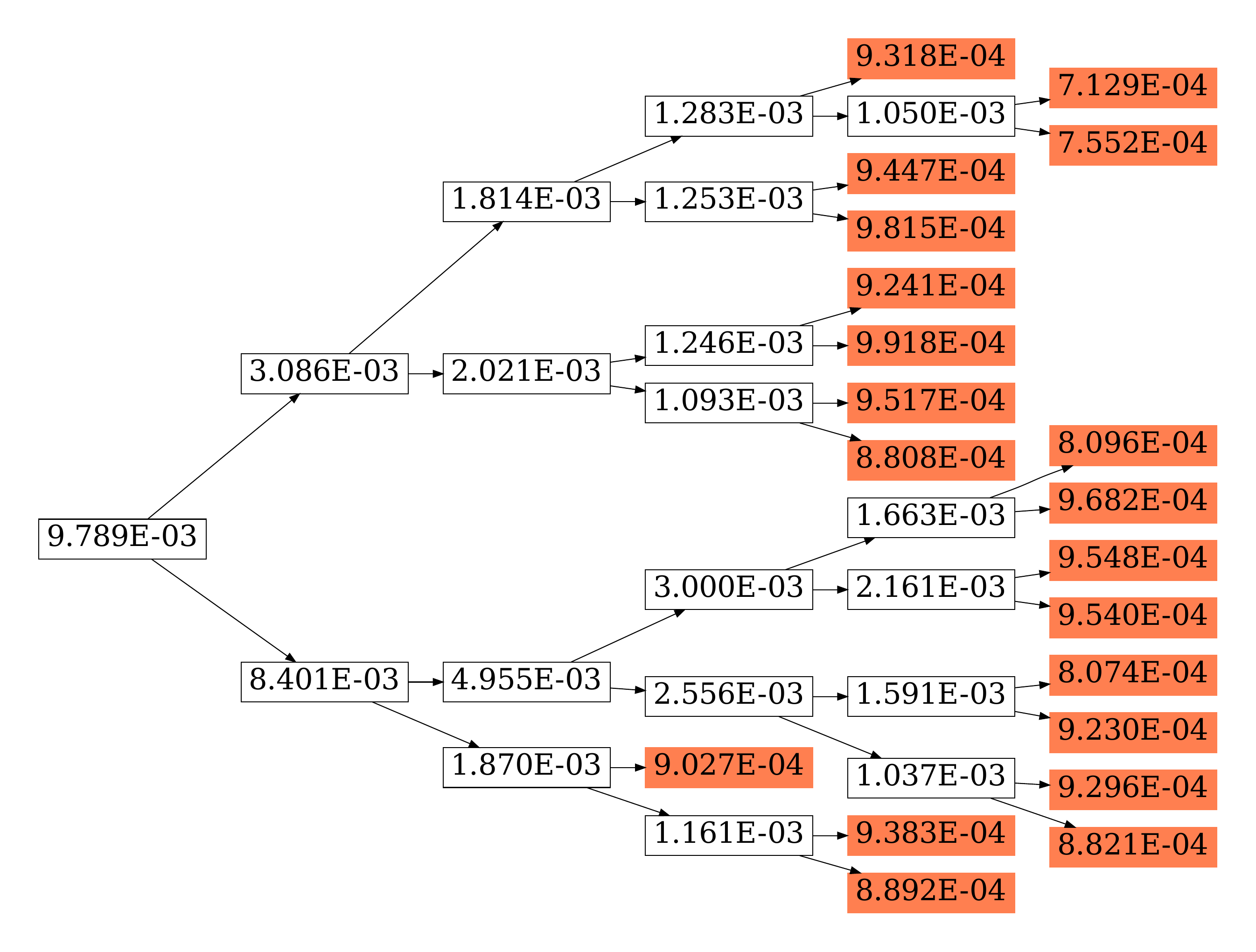}
\caption{\textbf{KdV problem:} Algorithm $\algoY$: accuracy $\dist(\cMalpha, \Valpha)$ for each node (leaves are colored in orange).}
\label{kdv:algoY-accuracy}
\end{figure}

\begin{figure*}[htbp]
\centering
\begin{subfigure}[t]{0.45\textwidth}
\centering
\includegraphics[width=\textwidth]{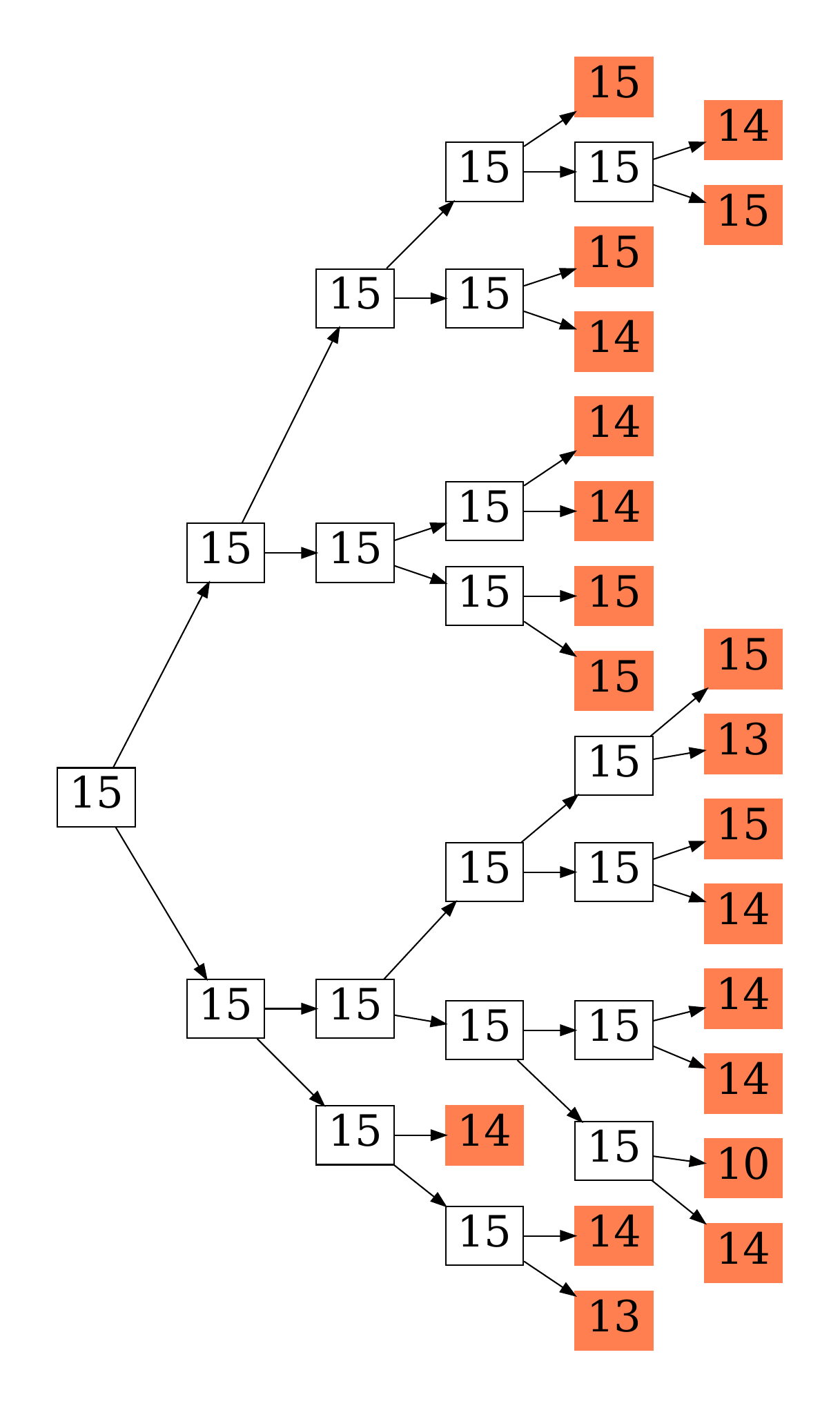}
\caption{Dimension $n(\alpha)=\dim(\Valpha)$ for each node.}
\label{kdv:algoY-dim}
\end{subfigure}%
~ 
\begin{subfigure}[t]{0.55\textwidth}
\centering
\includegraphics[width=\textwidth]{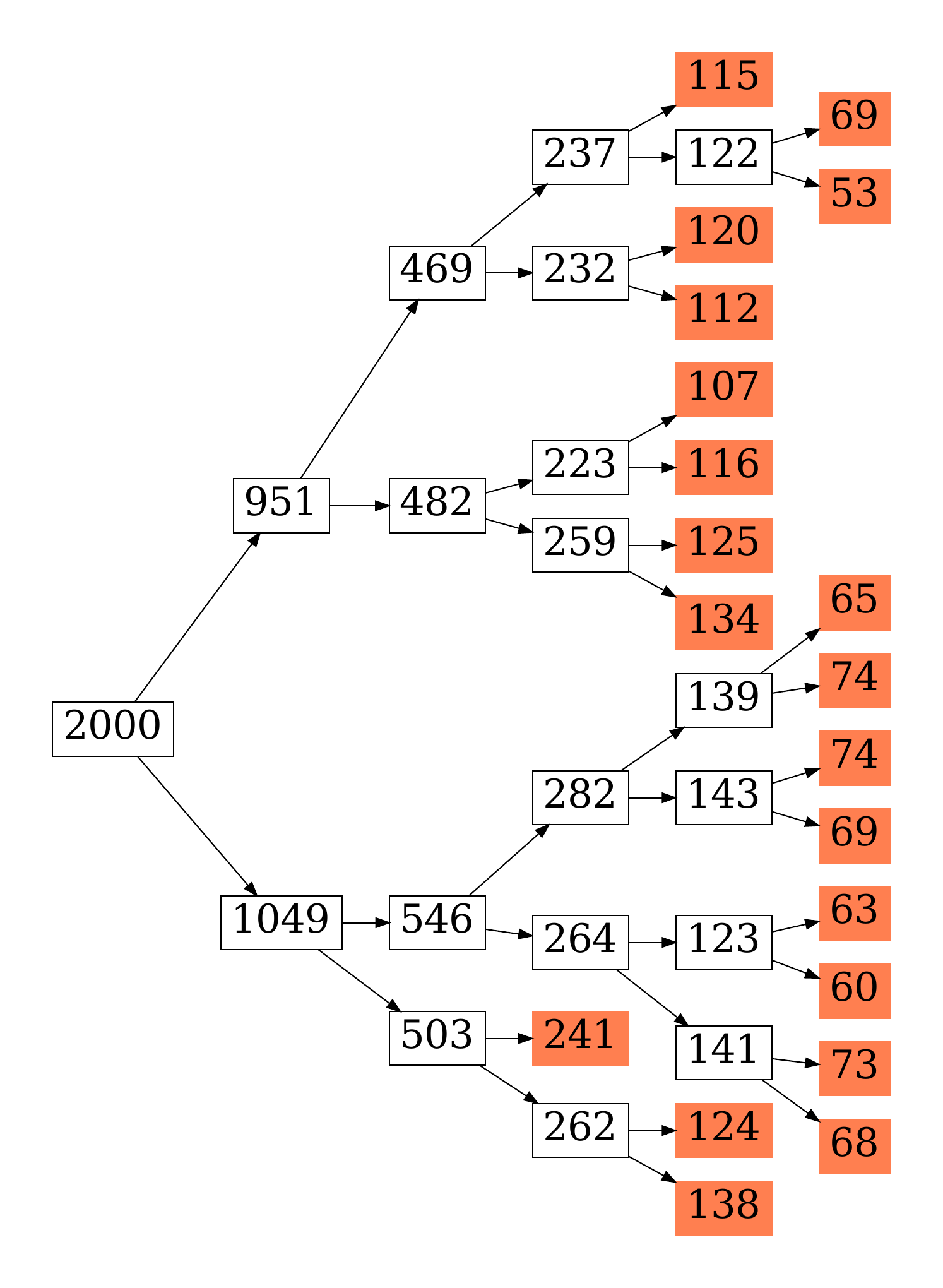}
\caption{Number of training snapshots for each node.}
\label{kdv:algoY-nb-snapshots}
\end{subfigure}
\caption{\textbf{KdV problem:} Algorithm $\algoY$: dimension and number of training snapshots in each tree node (leaves are colored in orange).}
\end{figure*}

To make a comparison with the \algoY~algorithm, we run it with $n_{\max}=15$ which is the average dimension that the spaces from \algoM~reach at convergence. In this case, \algoY~builds a library with $20$ spaces whose dimensions are $\{10,13,14,15\}$ and the amount of spaces with each dimension is $\{1,2,10,7\}$. The spaces of the library are generated with 281 snapshots. Interestingly, even if this library is smaller in terms of cardinality with respect to the one built by \algoM, it requires more snapshots to be created. The detailed behavior of the \algoY~algorithm is summarized in Figures \ref{kdv:algoY-accuracy}, \ref{kdv:algoY-dim} and \ref{kdv:algoY-nb-snapshots}.

\subsection{Main conclusion of the numerical tests}
\label{sec:numerics:conclusions}
The summary of Table \ref{numerics:maintable} illustrates that the superiority of one algorithm over another is relative, and depends on the criterion that is used. As predictited by the theory, the \algoPG~algorithm presents a fast decay for elliptic problems, and it is the method of choice when coercivity is not degenerate (see results for \diffone). As soon as one deviates from this setting, the other methods start becoming competitive. If one needs to work with local spaces of small dimension to reduce the online cost, \algoY~seems to be a good choice in most cases. This comes at the cost of having to store more snapshots to define the library (see the examples \difftwo, \convdiff~and \kdv). Thanks to its tree structure (namely one snaphots per node), the library generated by \algoM~requires less snapshots to generate the spaces but their dimension is often slightly larger than the ones obtained with \algoY~in our examples. As already brought up, \algoM~also comes with the advantage of being able to deal with very general parameter domains $\rY$ but the mapping $\iota$ of \eqref{eq:iota-algoM} is not trivial to deal with in practice.

\section{Conclusion}
The two different tree-based library approximations that we have explored in this work are very general, and can be applied to Hilbert, Banach and Wasserstein spaces. They can be formulated with minor modifications in general metric spaces. The most novel approach, called \algoM~in the paper, appears to be a good trade-off between the amount of snapshots that are needed to define the spaces of the library, and their respective dimensions. More broadly, the numerical experiments show a mixed landscape where each algorithm has relative advantages with respect to the others. A theoretical study on the convergence rate of both \algoY~and \algoM~is not entirely trivial and left for future work. Another element left for future work is the practical implementation of the mapping $\iota$ of \eqref{eq:iota-algoM} for $\algoM$.

\section*{Acknowledgments and disclosure of funding}
The authors would like to thank Andrea Bonito, Albert Cohen, Ron DeVore, Peter Jantsch and Guergana Petrova for fruitful discussions. \\
This research is part of the programme DesCartes and is supported by the National Research Foundation, Prime Minister's Office, Singapore under its Campus for Research Excellence and Technological Enterprise (CREATE) programme. It is also supported by the NSERC Grant RGPIN-2021-04311.

\bibliographystyle{unsrt}
\bibliography{references}

\end{document}